\documentclass{IEEEtran}
\usepackage[usenames,dvipsnames]{xcolor}

\usepackage{graphicx,subfig}
\usepackage{amsmath,amssymb,amsthm}
\usepackage{enumitem}

\def\conv{\mathrm{conv}}

\def\Prob{\mathrm{Prob}}
\def\1{\mathbf{1}}
\def\U{\mathcal{U}}

\def\P{\mathcal{P}}
\newcommand{\vv}[1]{\boldsymbol{#1}}
\newcommand{\vvi} [2]{ \boldsymbol{#1}^{(#2)} }
\newcommand{\vi} [2]{ {#1}^{(#2)} }

\newcommand{\vvu} [2]{ \boldsymbol{#1}^{[#2]} }
\newcommand{\dl}[1]{{#1}}

\newtheorem{definition}{Definition}
\newtheorem{corollary}{Corollary}

\newtheorem{theorem}{Theorem}
\newtheorem{lemma}{Lemma}

\theoremstyle{remark}

\usepackage{varwidth,calc}
\newlength\TableWidth%
\newsavebox\mybox%
\newcommand{\dynwidth}[2][(999)]{%
\begin{lrbox}{\mybox}%
\begin{varwidth}{10.0\columnwidth}
$\displaystyle  #2$
\end{varwidth}
\end{lrbox}%
\settowidth\TableWidth{\usebox\mybox}
\ifnum\TableWidth>\columnwidth\setlength\TableWidth{\columnwidth}\fi
\resizebox{1.0\TableWidth-\widthof{#1}}{!}{\mbox{$\displaystyle #2$}
}%
}

\title{Descent With Approximate Multipliers is Enough: Generalising Max-Weight}
\author{V\'ictor Valls, Douglas J. Leith\\Trinity College Dublin\thanks{This work
  was supported by Science Foundation Ireland
  under Grant No. 11/PI/1177.}}

\begin{document}
\maketitle

\begin{abstract} 
We study the use of approximate Lagrange multipliers and discrete actions in solving convex optimisation problems.   We observe that descent, which can  be ensured using a wide range of approaches (gradient, subgradient, Newton, \emph{etc}.), is orthogonal to the choice of multipliers.   Using the Skorokhod representation for a queueing process we show that approximate multipliers can be constructed in a number of ways.  These observations lead to the generalisation of (i) essentially any descent method to encompass use of discrete actions and queues and (ii) max-weight scheduling to encompass new descent methods including those with unsynchronised updates such as block coordinate descent.   This also allows consideration of communication delays and of updates at varying time-scales within the same clean and consistent framework.
\end{abstract}


\section{Introduction}

\IEEEPARstart{I}{n}  this paper we study the use of approximate Lagrange multipliers and discrete actions in solving convex optimisation problems.  One of the motivations of considering an approximation of the Lagrange multipliers is that in some problems the exact Lagrange multipliers might not be available, but a \emph{noisy}, \emph{delayed} or \emph{perturbed} version of the multipliers is available instead. For instance, in distributed optimisation a node might not have access to the exact Lagrange multiplier in the system due to transmission delays or losses.  Another example is in network problems where discrete valued queue occupancies can be identified with approximate scaled Lagrange multipliers. Restriction to a discrete set of actions is motivated by the fact that in many network problems the actual decision variables are discrete valued {e.g.}, do we transmit a packet or not.   

Our work is rooted in max-weight (also referred to as backpressure) scheduling  approaches for queueing networks, which have been the subject of much interest for solving utility optimisation problems in a distributed manner.   Appealing features of the max-weight scheduling approach include the lack of a requirement for \emph{a priori} knowledge of the packet arrival process and its support for a discrete action set. Our interest is not only in making these sorts of features available beyond queueing network applications but also in generalising max-weight approaches to allow advantage to be taken of the wealth of methods that exist for convex optimisation.


The main contributions are, in summary as follows.   
\begin{enumerate}[leftmargin=*,align=left,itemindent=.5cm,labelwidth=\itemindent]
\item Derivation of new conditions for ensuring convergence to a ball around the optimum of a convex optimisation, namely that descent with approximate multipliers is enough.    Descent is orthogonal to the choice of multipliers and can be ensured using a wide range of approaches (gradient, subgradient, Newton, \emph{etc}.).   Similarly, there exist a number of ways to construct approximate multipliers.
\item Regarding descent, we introduce a number of new convergence results for the unconstrained minimisation of a time-varying convex function when confined to selecting update steps from a specified action set, which may be discrete.  These extend the results in \cite{valls2014max} to allow unsynchronised updates, thereby generalising block co-ordinate approaches to encompass a broad class of discrete actions.   Unsynchronised updates are of interest for many reasons, not least being their greater efficiency in large-scale problems. 
\item Regarding approximate multipliers, using the Skorokhod representation for a queueing process the accuracy of approximate multipliers can be directly related to the difference between the exact and approximate queue increments as measured via the Skorokhod metric.  While max-weight uses a running average to control the approximation error,  we show that other methods can also be used.   This leads to an immediate generalisation of essentially any descent method to encompass use of discrete actions and queues. 
\end{enumerate}
Our analysis uses only elementary methods, requiring neither sophisticated Foster-Lyapunov nor fluid-limit machinery.

\subsection{Related Work}
Max-weight scheduling was introduced by Tassiulas and Ephremides in their seminal paper \cite{tassiulasstability}.  
Independently, \cite{stolyargreedy,eryilmazfair,neelypower} proposed extensions to the max-weight approach to accommodate concave utility functions.   
In \cite{stolyargreedy} the \emph{greedy primal-dual} algorithm is introduced for network linear constraints and  utility function which is continuously differentiable and concave. This is extended in \cite{stolyar2006greedy} to consider continuously differentiable nonlinear constraints.
In \cite{eryilmazfair} the utility fair allocation of throughput in a cellular downlink is considered.   
The work in \cite{neelypower} considers power allocation in a multi-beam downlink satellite communication link with the aim of maximising throughput while ensuring queue stability.  
This is extended in \cite{neelydynamic,neelyenergy, neelyfairness,neelybook} to develop the \textit{drift plus penalty} approach. 

With regard to the existence of a connection between the discrete-valued queue occupancy in a queueing network and continuous-valued Lagrange multipliers, this has been noted by several authors, see for example \cite{neely05,lin06}, but we are aware of few rigorous results before the recent work in \cite{valls2014max}.  Notable exceptions include \cite{neely11}, which establishes that a discrete queue update tends on average to drift towards the optimal multiplier value, and \cite{stolyargreedy} which shows for the \emph{greedy primal-dual} algorithm that asymptotically as design parameter $\beta \rightarrow 0$ and $t \rightarrow \infty$ the scaled queue occupancy converges to the set of dual optima.

Selection of a sequence of actions in a discrete-like manner is also considered in the convex optimisation literature.  The \emph{nonlinear Gauss-Seidel} algorithm, also known as \emph{block coordinate descent} \cite{bertsekasnonlinear,bertsekasparallel} minimises a convex function over a convex set by updating one co-ordinate at a time.  The convex function is required to be continuously differentiable and strictly convex and, unlike in the max-weight algorithms discussed above, the action set is convex.  The classical Frank-Wolfe algorithm \cite{frankwolfe} also minimises a convex continuously differentiable function over a polytope by selecting from a discrete set of descent directions, although a continuous-valued line search is usually used to determine the final update.   

\subsection{Notation}
Vectors and matrices are indicated in bold type.  Since we often use subscripts to indicate elements in a sequence, to avoid confusion we usually use a superscript $\vi{x}{i}$ to denote the $i$'th element of a vector $\vv{x}$.   With a mild abuse of notation where there is no scope for confusion we also sometimes use $\vi{f}{u}$, where $u$ is a set, to index a collection of functions {e.g.}, $\vi{f}{u}, u\in\U\subset2^{\{1,\cdots,n\}}$.
The $i$'th element of operator $[\vv{x}]^{[0,\bar{\lambda}]}$ equals $\vi{x}{i}$ (the $i$'th element of $\vv{x}$) when $\vi{x}{i}\in{[0,\bar{\lambda}]}$ and otherwise equals $0$ when $\vi{x}{i}<0$ and $\bar{\lambda}$ when $\vi{x}{i}\ge\bar{\lambda}$.  Note that we allow $\bar{\lambda}=+\infty$, and following standard notation in this case usually write $[x]^+$ instead of $[x]^{[0,\infty)}$.
 The subgradient of a convex function $f$ at point $\vv{x}$ is denoted $\partial f(\vv{x})$. {For two vectors $\vv{x},\vv{y} \in \mathbb R^m$ we use element-wise comparisons $\vv{x} \succeq \vv{y}$ and $\vv{y} \succ \vv{x}$ to denote when $\vi{y}{i} \ge \vi{x}{i}$, $\vi{y}{i} > \vi{x}{i}$ respectively for all $i=1,\dots,m$.}   
We let $\vvu{x}{u}\in \mathbb{R}^{|u|}$ denote the vector consisting of the subset  $u$ of elements of vector $\vv{x}\in\mathbb{R}^n$.   For example, if $u=\{2,3\}$ then $\vvu{x}{u}=[\vi{x}{2}, \vi{x}{3}]^T$.  
We also denote by $\vv{U}_{u}:\mathbb{R}^n\rightarrow \mathbb{R}^n$ the mapping such that vector $\vv{U}_{u}\vv{x}$ has elements in set $u$ equal to the corresponding elements in $\vv{x}$ and all other elements equal to $0$. That is, $(\vv{U}_{u}\vv{x})^{[u]}=\vvu{x}{u}$ and $\vi{(\vv{U}_{u}\vv{u})}{j}=0$, $j\in\{1,\cdots,n\}\setminus u$.  For example, if $u=\{2\}$ and $n=3$ then $\vv{U}_{u}\vv{x}=[0,\vi{x}{2}, 0]^T$.


\section{Preliminaries}

\subsection{Convex Optimisation}
Let $P$ denote the following convex optimisation problem:  
\begin{align*}
& \underset{\vv{z} \in C}{\text{minimise }} \qquad  f(\vv{z}) \\
& \underset{}{\text{subject to}} \qquad  \vv{g}(\vv{z}) \preceq \vv{0}
\end{align*}
where $f,\vi{g}{j}: \mathbb R^n \to \mathbb R$, $j=1,\dots,m$ are convex functions, $\vv{g}(\vv{z}) = [\vi{g}{1}(\vv{z}), \dots, \vi{g}{m}(\vv{z})]^T$ and $C$ a convex subset of $\mathbb R^n$. Let $C_0 := \{ \vv{z} \in C \mid \vv{g}(\vv{z}) \preceq \vv{0} \}$ and we will assume it is non-empty, {i.e.}, problem $P$ is feasible. Further, we will denote by $C^\star := \arg \min_{\vv{z} \in C_0} f(\vv{z}) \subseteq C_0$ the set of optima and $f^\star := f(\vv{z}^\star)$, $\vv{z}^\star \in C^\star$.

The Lagrange penalty function $L: \mathbb R^{n} \times \mathbb R^m \to \mathbb R$ associated with problem $P$ is 
\begin{align}
L(\vv{z},\vv{\lambda}) := f(\vv{z}) + \vv{\lambda}^T \vv{g}(\vv{z})
\end{align}where $\vv{\lambda} \in \mathbb R^m_+$.  Note that $L(\vv{z}, \vv{\lambda})$ is convex in $\vv{z}$ for a fixed $\vv{\lambda}$ and linear in $\vv{\lambda}$ for a fixed $\vv{z}$. The dual function is
\begin{align*}
{q}(\vv{\lambda}) := 
L(\vv{z}^*(\vv{\lambda}), \vv{\lambda}) 
= \min_{\vv{z} \in C} L(\vv{z}, \vv{\lambda}) && \vv{\lambda} \succeq \vv{0}
\end{align*}
where $\vv{z}^*(\vv{\lambda}) \in \arg \min_{\vv{z} \in C} L(\vv{z}, \vv{\lambda})$. Importantly, since $q(\vv{\lambda})$ is the minimum of a collection of affine functions it is concave and $q(\vv{\lambda}) \le L(\vv{z}, \vv{\lambda})$ for all $\vv{\lambda} \succeq \vv{0}$.   
The dual function is also Lipschitz continuous when $C$ is closed and bounded.  To see this observe that \begin{align}
{q}(\vv{\mu}) - {q}(\vv{\lambda})
&\le L(\vv{z}^*(\vv{\lambda}),\vv{\mu}) - L(\vv{z}^*(\vv{\lambda}),\vv{\lambda})\\
&= (\vv{\mu}-\vv{\lambda})^T\vv{g}(\vv{z}^*(\vv{\lambda}))\\
&\le \|\vv{\mu}-\vv{\lambda}\|_2 \| \vv{g}(\vv{z}^*(\vv{\lambda}))\|_2 
\le \|\vv{\mu}-\vv{\lambda}\|_\infty m{\bar g} \label{eq:qlip}
\end{align}
where ${\bar g}:= \max_{\vv{z} \in C} \|\vv{g}(\vv{z})\|_\infty$ \dl{is finite when $C$ is closed and bounded as $\vv{g}$ convex (so continuous)}.  Similarly, $ {q}(\vv{\lambda}) - {q}(\vv{\mu})\le  \|\vv{\mu}-\vv{\lambda}\|_\infty m{\bar g}$ and so $|{q}(\vv{\mu}) - {q}(\vv{\lambda})| \le \|\vv{\mu}-\vv{\lambda}\|_\infty m{\bar g}$.  That is, $q(\cdot)$ is Lipschitz continuous with Lipschitz constant $m{\bar g} $.   

The dual function is of particular interest when we have strong duality, {i.e.}, the solution of the dual problem 
\begin{align}
\underset{\vv{\lambda} \succeq \vv{0}}{\text{maximise}} \quad q(\vv{\lambda})\label{eq:dualp}
\end{align}
and the solution to problem $P$ coincide, $q(\vv{\lambda}^\star) = f^\star$ where  $\vv{\lambda}^\star \in\arg\max_{\vv{\lambda} \succeq \vv{0}} q(\vv{\lambda})$.   A sufficient condition for strong duality to hold is that the Slater condition is satisfied:
\begin{definition}[Slater condition] \label{th:slater} Set $C_0 := \{ \vv{z} \in C \mid \vv{g}(\vv{z}) \preceq \vv{0} \}$ has non-empty relative interior, {i.e.}, there exists a point $\bar{\vv{z}} \in C$ such that $\vv{g}(\bar{\vv{z}}) \prec \vv{0}$.
\end{definition}
%
 %
%
\noindent The following corresponds to Lemma 1 in \cite{nedicprimal} and is a direct consequence of the Slater condition. 
\begin{lemma}[Bounded level sets]\label{th:setq} Let $Q_{\delta} := \{\vv{\lambda} \succeq \vv{0}: q(\vv{\lambda}) \ge q(\vv{\lambda}^\star) - \delta \}$ with $\delta \ge 0$ and let the Slater condition hold, \emph{i.e.,} there exists a vector $\bar{\vv{z}} \in C$ such that $\vv{g}(\bar{\vv{z}}) \prec \vv{0}$. Then, for every $\vv{\lambda} \in Q_\delta $ we have that 
\begin{align}
\|\vv{\lambda} \|_2 \le \mathcal Q :=  \frac{1}{\upsilon} (f(\bar{\vv{z}}) - f^\star + \delta) \label{eq:setqbound}
\end{align}
where $\upsilon := \min_{j \in \{1,\dots,m\}} - \vi{g}{j}(\bar{\vv{z}})$. 
\end{lemma}
\noindent \dl{The importance of Lemma 1 is that it establishes that the set of dual optima is bounded when the Slater condition holds}.  \dl{It follows immediately that} 
\begin{align*}
\max_{\bar \lambda \vv{1} \succeq \vv{\lambda}\succeq \vv{0}} q(\vv{\lambda}) 
= \max_{\vv{\lambda}\succeq \vv{0}} \ q(\vv{\lambda}) 
= f^\star
\end{align*}
where $\bar \lambda \le m \mathcal Q$. 

 

\subsection{Optimisation Problems $P^\prime$ and $P^\dagger$}
In addition to optimisation problem $P$ we will also consider two special cases.  The first, which we denote $P^\prime$, is where:
\begin{enumerate}
\item $C:=\conv(D)$, $D$ a compact subset of $\mathbb{R}^n$.    
\end{enumerate}
\dl{The second, which we denote $P^\dagger$, assumes \emph{in addition} that:
\begin{enumerate}
\item Objective $f$ and constraints $\vi{g}{i}$, $i=1,\dots,m$ have bounded curvature on $C$ {{with constants $\mu_f$ and $\mu_{\vi{g}{i}}$, $i=1,\dots,m$.}}  \dl{When $\vv{\lambda}\preceq \bar \lambda\1$ it follows that the Lagrangian has bounded curvature on $C$ with constant $\mu_L =\mu_f +\bar{\lambda} \vv{1}^T \vv{\mu}_{\vv{g}}$ where $\vv{\mu}_{\vv{g}} = [\mu_{\vi{g}{1}},\dots,\mu_{\vi{g}{m}}]^T$}\footnote{In the special when $\vv{g}$ is linear then curvature constant $\vv{\mu}_g=\vv{0}$ and there is no need for $\vv{\lambda}\preceq \bar{\lambda}\1$ to ensure the Lagrangian has bounded curvature.}.
\item The objective and constraints are $\U$-separable, where $\U\subset 2^{\{1,\cdots,n\}}$.  That is, $f(\vv{z})=\sum_{u\in\U} \vi{f}{u}(\vvu{z}{u})$ and $\vv{g}(\vv{z}) = \sum_{u\in\U} \vvi{g}{u}(\vvu{z}{u})$ 
with $\vi{f}{u}:\mathbb{R}^{|u|}\rightarrow\mathbb{R}$, $\vvi{g}{u}:\mathbb{R}^{|u|}\rightarrow\mathbb{R}^m$.   
It follows that the Lagrangian is also $\U$-separable $L(\vv{z},\vv{\lambda})=\sum_{u\in\U} \vi{L}{u}(\vvu{z}{u},\vv{\lambda})$. 
\item Set $D$ is $\U$-feasible, see Definition \ref{def:d} below.

\end{enumerate}
}

\subsection{Bounded Curvature \& Descent}

Bounded curvature ensures that the level sets of a convex function have a smooth boundary and so it is relatively easy to find a descent direction.   
\begin{definition}[Bounded curvature]
Let $h:M\rightarrow\mathbb{R}$ be  a convex function defined on domain $M\subset\mathbb{R}^n$.   We say the $h(\cdot)$ has bounded curvature on set $C\subset M$ if for any points $\vv{z}, \vv{z}+\vv{\delta} \in C$
\begin{align}\label{eq:bound}
h(\vv{z}+\vv{\delta}) - h(\vv{z}) \le    \partial h(\vv{z})^T \vv{\delta} + \mu_h \| \vv{\delta} \|_2^2
\end{align}
where {$\mu_h \ge 0$} is a constant that does not depend on $\vv{z}$ or $\vv{\delta}$.  
\end{definition}
%
%
%
\noindent The following generalises Lemma 8 in \cite{valls2014max}.
\begin{lemma}[Discrete Descent]\label{lem:two}
Let $D$ be a compact subset of $\mathbb{R}^n$.  For any feasible point $\vv{y}\in C=\conv(D)$ and any vector $\vv{z}\in\mathbb{R}^n$ there exists a point $\vv{x}\in D$ such that $\vv{z}^T(\vv{x}-\vv{y}) \le 0$.
\end{lemma}
\begin{IEEEproof}
By Carath\'eŽodory's theorem, a point $\vv{y}\in C=\conv(D)\subset \mathbb{R}^n$ can be written as the convex combination of at most $n+1$ points $D^\prime:=\{\vv{x}_1,\dots,\vv{x}_{|D^\prime|}\}$ from $D$.  That is, $\vv{y}=\sum_{i=1}^{|D^\prime|}Ê\vi{a}{i} \vv{x}_i$ with $\sum_{i=1}^{|D^\prime|} \vi{a}{i}=1$, $\vi{a}{i}\in[0,1]$.  Hence, $\vv{z}^T(\vv{x}-\vv{y})=\sum_{i=1}^{|D^\prime|} \vi{a}{i}\vv{z}^T(\vv{x}-\vv{x}_i)$ for $\vv{x}\in D^\prime$.  Select $\vv{x} \in \arg \min_{\vv{w}\in D^\prime} \vv{z}^T\vv{w}$ (the min always exists since $D^\prime$ is finite).     Then $\vv{z}^T\vv{x} \le \vv{z}^T\vv{x}_i$ for all $\vv{x}_i\in D^\prime$ and so $\vv{z}^T(\vv{x}-\vv{y})\le 0$.
\end{IEEEproof}
Lemma \ref{lem:two} is fundamental and establishes that descent is possible even when the set of available directions is highly constrained {i.e.}, it is restricted to set $D$, which might be just the extreme points of $C$, a coordinate basis for $C$, \emph{etc}.   \dl{When $D$ consists of a finite set of points then Lemma \ref{lem:two} tells us that a descent direction can be found by direct search, without any need to compute the gradient.}


\subsection{Queue Continuity}\label{sec:cont}

Queue continuity plays a key role in our analysis of approximate multipliers.   


\dl{The following lemma, which corresponds to Proposition 3.1.2(i) in \cite{meyn2008control}, provides a useful way to represent a queueing process.} 
\begin{lemma}[Skorokhod map]\label{th:lindleys}
Consider the sequence $\lambda_{k+1} = [\lambda_k + \ x_k]^+$ for $k=1,2,\dots$ with $\lambda_1 
\ge 0$, $x_k \in \mathbb R$. Then,
\begin{align}
 \lambda_{k+1} 
&=  \max\Big\{ \max_{1 \le j \le k} \sum_{i=j}^k x_{i}, \Big[\sum_{i=1}^k x_i  + \lambda_1 \Big]^+  \Big\} \label{eq:fulllindleys}
\end{align}
When $\lambda_1=0$ this simplifies to $\lambda_{k+1} =  [ \max_{1 \le j \le k} \sum_{i=j}^kx_{i} ]^+ $.
\end{lemma}

\dl{Consider the sequences ${\lambda}_{k+1} = [{\lambda}_k + x_k]^+$ and ${\mu}_{k+1} = [{\mu}_k + y_k]^+$ with increments $x_k, y_k\in\mathbb{R}$ and $\lambda_1=\mu_1$, $k=1,2\dots$.    These updates can be thought of as fluid-like real-valued queues with real-valued increments.  Note that queues are more usually considered to be integer valued and we will return to this later.    The folllowing continuity result follows immediately from Lemma \ref{th:lindleys} and corresponds to Proposition 3.1.2(ii) in \cite{meyn2008control}.}
\begin{lemma}[Queue Continuity] \label{th:queuecontinuity}
Consider sequences $\{x_k\}$ and $\{ y_k\}$ of points from $\mathbb R$ and updates $\lambda_{k+1} = [\lambda_k +  x_k]^+$, $ \mu_{k+1} = [ \mu_k +  y_k]^+$ with $\lambda_1  = \mu_1\ge 0$. Then, 
\begin{align}
| \lambda_{k+1}  - \mu_{k+1}| 
&\le  2 \max_{1 \le j \le k }  \Big| \sum_{i=1}^j (x_i - y_i) \Big| \label{eq:qcont1}
\end{align}
\end{lemma}
%

\noindent  \dl{Hence, when $|\sum_{i=1}^k (x_i - y_i)| \le \sigma_0$ for all $k=1,2,\cdots$ then $| \lambda_{k}  - \mu_{k}|  \le 2\sigma_0$ for all $k=1,2,\cdots$.  Lemma \ref{th:queuecontinuity} states that the queue occupancy is Lipschitz continuous in the sequence of increments, with Lipschitz constant $2$, when the Skorokhod metric is used to measure the distance between the sequences of increments (namely $|\sum_{i=1}^k (x_i - y_i)|$). }
%
%
Note that in the special case where ${\lambda}_k = 0$ $\forall k$ then $| \lambda_{k+1}  - \mu_{k+1}| \stackrel{(a)}{=} \mu_{k+1} =  [\max_{1 \le j \le k} \sum_{i=j}^ky_{i}]^+  $ where $(a)$ follows because $\mu_{k+1}\ge 0$.  Hence, in this case when $\sum_{i=1}^j y_{i} \le\sigma_0$ for all $j=1,2,\dots$ then $| \lambda_{k}  - \mu_{k}|  \le \sigma_0$.   The requirement that $\sum_{i=1}^j y_{i} \le \sigma_0$ is much weaker than the requirement that  $|\sum_{i=1}^j (x_i - y_i)| \le \sigma_0$.  In particular, it allows $\sum_{i=1}^j y_{i}$ to be unbounded below and so, for example, the mean value of $y_{i}$ may be negative yet $\mu_k$ stays close to ${\lambda}_k = 0$. 

\dl{
The following is due to Theorem 4 in \cite{valls2014max}.
\begin{lemma}[Running Average] \label{th:auxiliarymultiplier1}
Let $\{ {x}_k \}$ be an arbitrary sequence of points from $\mathbb{R}$ and
 \begin{align}
 {z}_{k+1}  = (1-\beta) {z}_k + \beta {x}_k\label{eq:running}
\end{align}
with $\beta \in (0,1)$, ${z}_1 = x_1$.   Then 
\begin{align}
\textstyle |\sum_{i=1}^k {z}_k - {x}_k| \le \sigma_1 / \beta
\end{align}
where constant $\sigma_1\ge 0$.
\end{lemma}
\noindent Lemma \ref{th:auxiliarymultiplier1} states that when sequence $\{z_k\}$ is a running average of sequence $\{x_k\}$, then the sequences are close in the Skorokhod metric.  The following useful result is now immediate.
\begin{lemma} \label{th:auxiliarymultiplier}
Consider updates 
\begin{align}
{\vv{\lambda}}_{k+1}  & = [ {\vv{\lambda}}_k + \alpha (\vv{A}{\vv{z}}_{k+1} - \vv{b})]^+ \\
{\vv{\mu}}_{k+1}  & =  [{\vv{\mu}}_k + \alpha (\vv{A}\vv{x}_{k} - \vv{b}_k)]^+
\end{align}
 where step size $\alpha>0$, matrix $\vv{A}\in\mathbb{R}^{m\times n}$, $\{ \vv{x}_k \}$ is  an arbitrary sequence of points from $D\subset\mathbb{R}^n$ and
 \begin{align}
 \vv{z}_{k+1}  = (1-\beta) \vv{z}_k + \beta \vv{x}_k\label{eq:running}
\end{align}
with $\beta \in (0,1)$, $\vv{z}_1 \in C:=\conv{(D)}$. Further, suppose that $\{\vv{b}_k\}$ is a sequence of points in $\mathbb R^m$ such that $| \sum_{i=1}^k (\vi{b}{j}_i - \vi{b}{j}) | \le \sigma_2$  for all $j=1,\dots,m$, $k=1,2,\dots$. Then,
\begin{align}\label{eq:close}
\| {\vv{\mu}}_{k} - \vv{{\lambda}}_{k}\|_2 \le 2m \alpha (\sigma_1 / \beta + \sigma_2) ,\quad k=1,2,\dots
\end{align}
where ${\sigma_1} :=  2 \max_{\vv{z} \in C} \| \vv{A} \vv{z} \|_\infty$. 
\end{lemma}
}
%
\noindent 
Note that while the updates in Lemma \ref{th:auxiliarymultiplier} involve scaling factor $\alpha$, by rescaling we can remove this factor.  Namely, letting $\vv{Q}_k=\vv{\mu}_k/\alpha$ then ${\vv{Q}}_{k+1}   = [ {\vv{Q}}_k + (\vv{A}{\vv{x}}_{k} - \vv{b})]^+$.


\begin{figure}
\centering
\includegraphics[width=0.6\columnwidth]{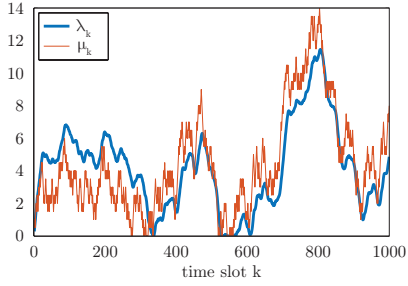}
\caption{Example realisations of $\lambda_k$ and $\mu_k$}\label{fig:ex0}
\end{figure}
Figure \ref{fig:ex0} shows an example of sequences $\lambda_k$ and $\mu_k$, $k=1,2,\dots$ where $A=1$, $b_k=b=0.5$, $\alpha=1$, $\beta=0.1$ and $x_k$ takes values independently and uniformly at random from set $\{0,1\}$.   It can be seen that the distance between $\lambda_k$ and $\mu_k$ remains uniformly bounded over time.

\section{Descent}
Given a sequence $F_k$, $k=1,2,\cdots$ of convex functions, in this section we consider how to construct a corresponding sequence of values $\vv{z}_k$ such that $F_k(\vv{z}_k)$ converges to an $\epsilon$-ball around its minimum.  That is, such that $F_k(\vv{z}_k) - \min_{\vv{z}\in C} F_k(\vv{z}) \le \epsilon$.   In particular, our interest is in sequences $\vv{z}_k$, $k=1,2,\cdots$ which are (i) constrained to be the running average (\ref{eq:running}) of an associated sequence $\vv{x}_k$, $k=1,2,\dots$ and where (ii) each $\vv{x}_k$ must be selected from a specified set $D$ of actions.  This will provide a basic building block for our later analysis.   

In addition to these constraints on the admissible choices of $\vv{z}_k$ and $\vv{x}_k$, we will also allow constraints on the set of elements of $\vv{z}_k$ that can be updated at each step $k$.   This allows us to not only capture lack of synchronism in the updates, {e.g.}, in a distributed setting, but also computational constraints, {e.g.}, where calculation of a limited update is much more efficient than that of a full update.

\subsection{Unsynchronised Updates}
 We introduce the following definition related to unsynchronised updates.
\begin{definition}[Admissible Update Set]\label{def:u}
Let $\bar \U\subset 2^{\{1,\dots,n\}}$ be a partition of $\{1,\dots,n\}$ {i.e.}, $\cup_{u\in\bar{\U}}u=\{1,\dots,n\}$, $u\cap v=\emptyset$ for $u$, $v\in\bar{\U}$ and $u\ne v$.  
Then we say $\U\subset 2^{\{1,\dots,n\}}$ is an \emph{admissible update set} if $\U=\bar{\mathcal U}$ or $\U=\bar {\mathcal U} \cup \{1,\dots,n\}$.
\end{definition}
\noindent That is $\U\subset 2^{\{1,\cdots,n\}}$ is an admissible update set if its elements are either the full set $\{1,\dots,n\}$ or belong to a partition of this full set.  Letting set $\U\subset 2^{\{1,\dots,n\}}$ be an admissible update set, let each $u\in \U$ define a set of elements $\vvu{z}{u}_k$ of vector $\vv{z}_k$ that can be updated jointly at time step $k$.   The technical conditions in Definition \ref{def:u} ensure that these updates are suitably well behaved.  Recalling that $\vv{U}_{u}\vv{x}$ is the vector with elements at the positions specified by $u$ equal to the corresponding elements in $\vv{x}$ and all other elements zero, we also introduce a second definition: 

\begin{definition}[$\U$-Feasible Domain]\label{def:d}
Let $D\subset\mathbb{R}^n$ and $C=\conv(D)$.   Let $B:=\{\vv{z}+ \vv{U}_{u} (\vv{x}-\vv{z}): \vv{z}\in C, \vv{x}\in D, u\in \U\subset 2^{\{1,\dots,n\}} \}$.  Then, set $D$ is said to be \emph{$\U$-feasible} if $B\subset C$.
\end{definition} 

\noindent  Letting $\vv{y}= \vv{z}+ \vv{U}_{u} (\vv{x}-\vv{z})$, $\vv{z} \in C$, $\vv{x}\in D$ then observe that $\vv{z}+ \beta \vv{U}_{u} (\vv{x}-\vv{z}) =(1-\beta) \vv{z} +\beta\vv{y}  \in C$.   That is, when set $D$ is $\U$-feasible then the running average (\ref{eq:running}) lies in set $C$ even when updates to $\vv{z}$ are confined to a subset $u$ of elements of $\vv{z}$.
%
%
It also implies that set $C$ has the product form $C=\Pi_{u\in\U}C_{u}$ with $C_u\subset \mathbb{R}^{|u|}$, $u\in\U$.

\emph{Example}.  The special case where only a single element of $\vv{z}$ can be updated at a time corresponds to selecting $\U=\{\{1\},\{2\},\dots,\{n\}\}$ and set $D=\{\vv{e}_1,\dots,\vv{e}_n\}\cup\{\1\}$, where $\vv{e}_i$ denotes the unit vector with all elements zero apart from element $i$ and $\1$ denotes the all ones vector.   Then set $C$ is the unit hypercube in $\mathbb{R}^n$.  Vector $\vv{z}+\vv{U}_{u}(\vv{x}-\vv{z})$, $u\in\U$ has all elements the same as vector $\vv{z}$ apart from element $i$ which has value either $0$ or $1$ depending on the choice of vector $\vv{x}\in D$.   For both of these values the vector lies in $C$ and so $D$ is $\U$-feasible. 

\begin{definition}[$\U$-Separable]
We say a function $h:\mathbb{R}^n\rightarrow\mathbb{R}$ is $\U$-separable when $h(\vv{z})=\sum_{u\in\U} \vi{h}{u}(\vvu{z}{u})$ with $\vi{h}{u}:\mathbb{R}^{|u|}\rightarrow\mathbb{R}$ and $\U\subset 2^{\{1,\dots,n\}}$.   
\end{definition}


\subsection{Descent Using An Action Set}

We have the following corollary to Lemma \ref{lem:two} for a set $u$ of elements:
\begin{corollary}\label{cor:one}
For any $\vv{y}\in C=\conv(D)$, $\vv{z}\in\mathbb{R}^{n}$ and $u\subseteq\{1,\dots,n\}$ there exists a point $\vv{x}\in D$, a compact subset of $\mathbb{R}^n$, such that $(\vvu{z}{u})^T(\vvu{x}{u}-\vvu{y}{u}) \le 0$.
\end{corollary}
\begin{IEEEproof}
Observe that $\vvu{y}{u}=\sum_{j} \vi{a}{j} \vvu{x}{u}_j$ for any $\vv{y}\in C$.
\end{IEEEproof}
\noindent Using this, together with bounded curvature, it follows that a choice of $\vv{x}\in D$ ensuring descent always exists:
\begin{lemma}[Descent]\label{th:coord}
Let $F:\mathbb{R}^n\rightarrow\mathbb{R}$ be a convex function with bounded curvature on $C$ with curvature constant $\mu_F$.  Suppose points $\vv{y}$, $\vv{z}\in C=\conv(D)$ exist such that $F(\vv{z}+\vv{U}_{u}(\vv{y}-\vv{z}))\le F(\vv{z})-\epsilon$, where $D$ is a $\U$-feasible, compact subset of $\mathbb{R}^n$, $\epsilon>0$, $u\in\mathcal{U}\subset 2^{\{1,\dots,n\}}$.  
Selecting
\begin{align}
\vv{x} &\in\arg\min_{\vv{w}\in D} F(\vv{z}+ \beta \vv{U}_{u}(\vv{w}-\vv{z})) \label{eq:d1}
\end{align}
then
\begin{align*}
F( \vv{z} + \beta \vv{U}_{u}(\vv{x}-\vv{z}))<F(\vv{z})-\gamma\beta\epsilon
\end{align*}
 provided $\beta\le (1-\gamma)\min\{\epsilon / (\mu_F\bar{x}^2_D),1\}$ $0<\gamma<1$ where $\bar{x}_D:=\max_{\vv{y},\vv{z}\in \conv(D)}\|\vv{y}-\vv{z}\|_2$.
\end{lemma}
\begin{IEEEproof}
See the appendix.
\end{IEEEproof}

Note that Lemma \ref{th:coord} requires set $D$ to be $\U$-feasible and compact but places no other restriction on $D$.  For example, we can select $D=C$ or we can select $D$ such that it contains just the extreme points of $C$.   \dl{When set $D$ is small and finite then the solution to (\ref{eq:d1}) can be found by direct search without requiring computation of the gradient of function $F$.   When $D=C$ then (\ref{eq:d1}) can be solved using any convenient convex optimisation method.}

Using Lemma \ref{th:coord} we can now show that, even when $\vv{z}_k$ is constrained to be a running average of actions $\vv{x}_k\in D$ and at each time step only a set of elements $u\in\mathcal{U}$ can be updated, it is nevertheless still possible to construct a sequence $\vv{z}_k$ that forces descent.

\begin{theorem}[Convergence]\label{th:unsynch}
Let $u_k$, $k=1,2,\dots$ be a sequence of elements from admissible update set $\mathcal{U}$.   Suppose we can partition these this sequence into disjoint intervals $\{k_i,\dots,k_{i+1}\}$, $i=1,2,\dots$ with $k_i> k_{i-1}$
and such that over each interval $\cup_{k\in\{k_j,\dots,k_{j+1}\}} u_k = \{1,\dots,n\}$.  Suppose the largest size of interval $N=\max_{j}(k_{j+1}-k_j+1)$ is finite.   

Let action set $D$ be $\U$-feasible, compact subset of $\mathbb R^n$.   Let $F_k$, $k=1,2,\dots$ be a sequence of $\U$-separable convex functions, $F_k(\vv{z})=\sum_{u\in\U} \vi{F_k}{u}(\vvu{z}{u})$ where $\vi{F_k}{u}:\mathbb{R}^{|u|}\rightarrow \mathbb{R}$ is convex with bounded curvature and curvature constant $\mu_F$.    Let $\{\vv{z}_k\}$ be a sequence of vectors satisfying
\begin{align}
\vv{x}_k &\in\arg\min_{\vv{x}\in D} \vi{F}{u_k}_k(\vvu{z}{u_k}_k+ \beta (\vvu{x}{u_k}-\vvu{z}{u_k}_k))\label{eq:update1}\\
\vv{z}_{k+1}&=\vv{z}_k+ \beta \vv{U}_{u_k} (\vv{x}_k-\vv{z}_k)
\end{align}
with $\vv{z}_1\in C$.

Suppose that parameter $\beta$ is sufficiently small that
\begin{align}
\beta\le (1-\gamma)\gamma\min\Big\{\frac{\epsilon^\prime}{N\mu_F\bar{x}^2_D},1\Big\} \label{eq:beta}
\end{align}
with $0<\gamma<1$, $\epsilon^\prime>0$, $\bar{x}_D:=\max_{\vv{y},\vv{z}\in \conv(D)}\|\vv{y}-\vv{z}\|_2$ and that functions $F_k$ change sufficiently slowly that
\begin{align*}
|F_{k+1}(\vv{z}) - F_k(\vv{z})| \le \frac{\gamma_1\gamma}{2N}\beta\epsilon^\prime,\ \forall \vv{z}\in C
\end{align*}
where $\gamma_1\in(0,\gamma/2)$.  
Then for $k_i$ sufficiently large we have that, 
\begin{align*}
F_{k_{i}}(\vv{z}_{k_i})-F_{k_{i}}(\vv{y}_{k_{i}}^*)\le \epsilon:=\left({|\U|}+ 1\right)\epsilon^\prime
\end{align*}
where $\vv{y}^*_k\in\arg\min_{\vv{z}\in C} F_k(\vv{z})$.
\end{theorem}
\begin{IEEEproof}
See the appendix.
\end{IEEEproof}

Note that Theorem \ref{th:unsynch} adds the restriction that $F_k$ is $\U$-separable. 
Although this restriction might be weakened somewhat, it does seem that some form of separability is essential to ensure that descent happens when only a subset of elements of $\vv{z}_k$ are updated at each step.   When all elements can be updated simultaneously then $\mathcal{U}$ contains a single element $\{1,\dots,n\}$ and there is no need for $F_k$ to be separable.

One of the main difficulties in establishing Theorem \ref{th:unsynch} is that at each step we are constrained to select an action $\vv{x}_k$ from set $D$ when updating the running average $\vv{z}_k$.   This means that there is no ``null'' action leaving $\vv{z}_k$ unchanged since that would require selecting $\vvu{x}{u_k}_k$ equal to $\vvu{z}{u_{k}}_k$ which is generally not possible (for example, $\vv{x}_k$ may be constrained to be discrete valued yet $\vv{z}_k$ is real-valued\footnote{\dl{In the special case where $D=C$ then of course this issue does not arise.}}).  Consequently, we sometimes may be forced to allow $F_k$ to increase.  Nevertheless, Theorem \ref{th:unsynch} shows it is always possible to select $\vv{z}_k$ such that at the sub-sequence of times $k_i$, $i=1,2,\dots$ the function $F_k$ converges into a $\epsilon$-ball around its minimum.

Update (\ref{eq:update1}) is not the only possible choice for ensuring descent.  For example, the following theorem gives an alternative choice:
\begin{theorem}[FW Convergence]\label{th:FWunsynch}
Suppose the conditions of Theorem \ref{th:unsynch} are satisfied and, in addition, that set $D$ contains the extreme points of set $C$.  Let $\{\vv{z}_k\}$ be a sequence of vectors satisfying
\begin{align}
\vv{x}_k &\in\arg\min_{\vv{x}\in D} \partial F(\vv{z}_k)^T\vv{U}_{u_{k}}\vv{x}\label{eq:updateFW}\\
\vv{z}_{k+1}&=\vv{z}_k+ \beta \vv{U}_{u_{k}} (\vv{x}_k-\vv{z}_k)
\end{align}
with $\vv{z}_1\in C$.
Then for $k_i$ sufficiently large we have that, 
\begin{align*}
F_{k_{i}}(\vv{z}_{k_i})-F_{k_{i}}(\vv{y}_{k_{i}}^*)\le \epsilon:= \left({|\U|}+ 1\right)\epsilon^\prime
\end{align*}
where $\vv{y}^*_k\in\arg\min_{\vv{z}\in C} F_k(\vv{z})$.
\end{theorem}
\begin{IEEEproof}
See the appendix.
\end{IEEEproof}

\noindent Update (\ref{eq:updateFW}) generalises the classical Frank-Wolfe update \cite{frankwolfe} to minimisation of a sequence of functions and to allow partial updating at each time step.

\section{Descent With Approximate Multipliers is Enough}

The previous section considers finding specific types of sequence that converge to an $\epsilon$-ball around an unconstrained minimum.  We now extend consideration to include constraints.  The following approximation result is the key to this extension.

\begin{theorem}[Approximate Solution] \label{th:maintheorem} Consider the convex optimisation $P$ and update
\begin{align}
\vv{\lambda}_{k+1} = [ \vv{\lambda}_k  + \alpha \vv{g}(\vv{z}_k)]^+ \label{eq:dualupdatetheorem}
\end{align}
with $\alpha > 0$.  Suppose the Slater condition is satisfied and also the following two conditions:
\begin{itemize}
\item[(i)] Descent. $\vv{z}_k \in  C_\epsilon (\vv{\mu}_k)$ where $C_\epsilon (\vv{\mu}_k) := \left\{ \vv{z} \in C \mid L(\vv{z},\vv{\mu}_k) - q(\vv{\mu}_k) \le \epsilon \right\}$ with $\epsilon \ge 0$; 
\item[(ii)] Approximate Multipliers.  $\| \vv{\lambda}_k - \vv{\mu}_k \|_\infty \le \alpha \sigma_0$ for all $k=1,2,\dots$ with $\sigma_0 \ge 0$.
\end{itemize}
Then, 
\begin{align}
 &- \frac{2 m \bar \lambda^2 }{\alpha k} -  {\alpha} m (\bar g^2 / 2 + \sigma_0 (1+ \bar g))  -\epsilon  \notag\\ 
 &\qquad\qquad\le f(\vv{z}^\diamond_k) -  f^\star 
 \le  \epsilon + \alpha m ({\bar g}^2 + \sigma_0(1+\bar g))  + \frac{ 3 m \bar \lambda^2}{2\alpha k}  \notag
\end{align}
where $ \vv{z}^\diamond_k  := k^{-1} \sum_{i=1}^k \vv{z}_i$, ${\bar g}:= \max_{\vv{z} \in C} \|\vv{g}(\vv{z})\|_\infty$ and $\bar \lambda$ is a positive constant that does not depend on $k$.
\end{theorem}
\begin{IEEEproof}
See the appendix.
\end{IEEEproof}
Note that it is the average, $\vv{z}^\diamond_k$, which converges to a ball around the optimum, rather than $\vv{z}_k$ itself.  By selecting $\alpha$ and $\epsilon$ sufficiently small, we can make this ball as small as we like as $k\rightarrow\infty$.  Observe also that the bounds in Theorem \ref{th:maintheorem} are not asymptotic but rather can be applied at finite times.   As $k\rightarrow\infty$ the bounds simplify to $-{\alpha} m (\bar g^2 / 2 + \sigma_0 (1+ \bar g))  -\epsilon   \le f(\vv{z}^\diamond_k) -  f^\star  \le  \epsilon + \alpha m ({\bar g}^2 + \sigma_0(1+\bar g))$.

\begin{corollary}
Consider the setup of Theorem \ref{th:maintheorem} where update (\ref{eq:dualupdatetheorem}) is replaced with 
\begin{align}
\vv{\lambda}_{k+1} = [\vv{\lambda}_k + \alpha \vv{g}(\vv{z}_k)]^{[0,\bar \lambda]}.
\end{align} 
and $\vv{\lambda}_1\preceq\bar{\lambda}\1$, where the $j$'th element of operator $[\vv{x}]^{[0,\bar{\lambda}]}$ equals $\vi{x}{j}$ when $\vi{x}{j}\in{[0,\bar{\lambda}]}$ and otherwise equals $0$ when $\vi{x}{j}<0$ and $\bar{\lambda}$ when $\vi{x}{j}\ge\bar{\lambda}$.   Then, the bound in Theorem \ref{th:maintheorem} holds provided $\bar \lambda \le m \mathcal Q$ where $\mathcal Q$ is given in Lemma \ref{th:setq}.
\end{corollary}
\begin{IEEEproof}
The result follows directly from the fact that the set of dual optima is bounded (see Lemma \ref{th:setq}) and that the difference $\| \vv{\lambda}_k - \vv{\lambda}^\star \|_2^2$ decreases monotonically until $\vv{\lambda}_k$ converges to a ball around $\vv{\lambda}^\star$ (see Lemma \ref{th:bmulti}).
\end{IEEEproof}

Theorem \ref{th:maintheorem} establishes convergence to an optimum provided conditions (i) and (ii) are satisfied.  Condition (i) is that sequence $\vv{z}_k$, $k=1,2,\dots$ ensures that $L(\vv{z}_k,\vv{\mu}_k)$ stays within an $\epsilon$-ball of its minimum.   This is essentially a descent condition {i.e.}, whenever $L(\vv{z}_k,\vv{\mu}_k)$ is outside the $\epsilon$-ball we must select $\vv{z}_k$ to decrease $L$ so as to bring it back within the ball.    Note that approximate multiplier $\vv{\mu}_k$ is used in $L(\vv{z}_k,\vv{\mu}_k)$, so there is no need to know the exact multiplier $\vv{\lambda}_k$ in order to determine a suitable $\vv{z}_k$.      Condition (ii) is that approximate multiplier $\vv{\mu}_k$ stays close to the exact multiplier $\vv{\lambda}_k$.   

Our convergence analysis of algorithms for constrained optimisation now reduces to establishing whether these two conditions, which we refer to as descent with approximate multipliers, are satisfied.

\subsection{Some Examples}

\subsubsection{Classical Dual Subgradient Approach}
This uses the following update to solve convex optimisation $P$:
\begin{align}
&\vv{z}_{k} \in \arg\min_{\vv{z}\in C} L(\vv{z},\vv{\lambda}_k)\\
&\vv{\lambda}_{k+1} = [ \vv{\lambda}_k  + \alpha \vv{g}(\vv{z}_k)]^+
\end{align}
This choice of $\vv{z}_k$ ensures $\vv{z}_k \in  C_\epsilon (\vv{\lambda}_k)$ with $\epsilon=0$ and so satisfies Theorem \ref{th:maintheorem} condition (i), while the choice of $\vv{\mu}_k=\vv{\lambda}_k$ trivially satisfies condition (ii).

Theorem \ref{th:maintheorem} allows the dual subgradient approach to be immediately extended to use unsynchronised updates.  To see this, suppose in problem $P$ that the Lagrangian is $\U$-separable, $L(\vv{z},\vv{\lambda})=\sum_{u\in\U} \vi{L}{u}(\vvu{z}{u},\vv{\lambda})$ where $\U\subset 2^{\{1,\dots,n\}}$ is an admissible update set, and set $C$ is $\U$-feasible.   Hence, $\min_{\vv{z}\in C} L(\vv{z},\vv{\lambda}_k)= \sum_{u\in\U} \min_{\vv{w}\in C_u} \vi{L}{u}(\vv{w},\vv{\lambda}_k)$ where $C= \prod_{u\in\U}C_u$.   Allowing each $\vi{L}{u}$ to use a different approximate multiplier $\vvi{\mu}{u}_k$ leads to the following update 
\begin{align}
&\vvu{z}{u}_{k} \in \underset{\vv{w} \in C_u}{\arg\min} \ \vi{L}{u}(\vv{w},\vvi{\mu}{u}_k) &&\ u\in\U \label{eq:rr1}\\
&\vv{\lambda}_{k+1} = [ \vv{\lambda}_k  + \alpha \vv{g}(\vv{z}_k)]^+ \label{eq:rr2}
\end{align}
%
Consider now the situation where at each step $k$ only a subset of the elements $\vvu{z}{u}_k$ are updated {i.e.}, updates are unsynchronised.   The elements $u$ which are updated use multiplier $\vvi{\mu}{u}_k=\vv{\lambda}_{k}$.  Those elements $\vvu{z}{v}_k$, $v\ne u$ which are not updated at step $k$ can be formally thought of as being updated using the old multiplier value $\vvi{\mu}{v}_k=\vv{\lambda}_{k-\vi{\tau}{v}_k}$ where $k-\vi{\tau}{v}_k$ is the time step where $\vvu{z}{v}_k$ was last updated.   Provided $\vi{\tau}{v}_k$ is uniformly upper bounded by $\sigma_1$ then $\|\vv{\lambda}_k - \vvi{\mu}{v}_k \|_\infty$ is unformly bounded by $\alpha \sigma_1\bar{g}$ and so  by Theorem \ref{th:maintheorem} updates (\ref{eq:rr1})-(\ref{eq:rr2}) converge to a ball around the optimum.  Note that there is no need for objective $f$ or constraints $\vv{g}$ to have bounded curvature.
  
%


\subsubsection{Frank-Wolfe Dual Subgradient Approach}
Consider the following update to solve convex optimisation $P^\prime$:
\begin{align}
&\vv{y}_{k} \in \arg\min_{\vv{z}\in C} \partial_{\vv{z}}L(\vv{z}_k,\vv{\lambda}_k)^T\vv{z}\label{eq:fw}\\
&\vv{z}_{k+1}=(1-\beta)\vv{z}_k+\beta \vv{y}_k\\
&\vv{\lambda}_{k+1} = [ \vv{\lambda}_k  + \alpha \vv{g}(\vv{z}_{k+1})]^{[0,\bar \lambda]}\label{eq:fw2}
\end{align}
%
Suppose that  we select step size $\alpha$ sufficiently small\footnote{This is always possible since $|L(\vv{z},\vv{\lambda}_{k+1})-L(\vv{z},\vv{\lambda}_k)| \le \|\vv{\lambda}_{k+1}-\vv{\lambda}_k\|_2 \bar{g} \le \alpha m \bar{g}^2$ where $\bar{g}:=\max_{\vv{z}\in C}\|\vv{g}(\vv{z})\|_2$.   Since $C$ is compact and $\vv{g}$ convex (so continuous) then $\bar{g}$ exists and is bounded.  Hence, selecting $\alpha\le\frac{\gamma_1\gamma}{2\bar{g}}\beta\epsilon^\prime$ is sufficient. } that
\begin{align}
|L(\vv{z},\vv{\lambda}_{k+1})-L(\vv{z},\vv{\lambda}_k)|&\le \frac{\gamma_1\gamma}{2}\beta\epsilon^\prime,\ \forall \vv{z}\in C\label{eq:slow}
\end{align}
Then the Lagrangian $L(\cdot,\vv{\lambda}_k)$ satisfies the conditions on $F_k$ in Theorem \ref{th:FWunsynch}.  Hence, provided step size $\beta$ is sufficiently small that it satisfies condition (\ref{eq:beta}), applying Theorem \ref{th:FWunsynch} with update set $\U$ containing the single element $u_1=\{1,\cdots,n\}$ then $\vv{z}_k\in C_\epsilon(\vv{\lambda}_k)$ for $k$ sufficiently large and so satisfies Theorem \ref{th:maintheorem} condition (i).   The choice of $\vv{\mu}_k=\vv{\lambda}_k$ trivially satisfies condition (ii).

\subsubsection{Classical Max-Weight}
Consider the following update to solve convex optimisation $P^\prime$ with linear constraints ($\vv{g}(\vv{z})=\vv{A}\vv{z}-\vv{b}$):
\begin{align}
&\vv{x}_k\in\arg\min_{\vv{x}\in D} \partial_{\vv{z}}L(\vv{z}_k,\vv{\mu}_k)^T \vv{x} \label{eq:mw}\\
&\vv{z}_{k+1}=(1-\beta)\vv{z}_k+\beta \vv{x}_k\\
&\vv{\mu}_{k+1} = [ \vv{\mu}_k  + \alpha (\vv{A}\vv{x}_k-\vv{b}_k)]^+\label{eq:mw2}
\end{align}
where $\|\sum_{i=1}^k\vv{b}_i -\vv{b}\|_\infty \le \sigma_2$ for all $k=1,2,\dots$.  

Note that update (\ref{eq:mw})-(\ref{eq:mw2}) is not written in the usual max-weight notation.  However, recall that $\partial_{\vv{z}}L(\vv{z}_k,\vv{\mu}_k)^T = \partial f(\vv{z}_k)+\vv{\mu}_k^T\vv{A}$.   Also that by defining $\vv{Q}_{k}=\vv{\mu}_k/\alpha$ then $\vv{Q}_{k+1} = [ \vv{Q}_k  + \vv{A}\vv{x}_k-\vv{b}_k]^+$.   Identifying cost $-f(\vv{z})$ with utility $U(\vv{z})$ which is to be maximised, it then follows that (\ref{eq:mw})-(\ref{eq:mw2}) is can be rewritten equivalently as
\begin{align}
&\vv{x}_k\in\arg\max_{\vv{x}\in D} (\partial U(\vv{z}_k)-\alpha\vv{Q}_k^T\vv{A}) \vv{x} \\
&\vv{z}_{k+1}=(1-\beta)\vv{z}_k+\beta \vv{x}_k\\
&\vv{Q}_{k+1} = [ \vv{Q}_k  + \vv{A}\vv{x}_k-\vv{b}_k]^+
\end{align}
which is precisely the max-weight update of Stolyar \cite{stolyargreedy}.

The reason to write the max-weight update in the form (\ref{eq:mw})-(\ref{eq:mw2}) is to highlight its close connection with other convex optimisation approaches.  In particular, observe the close relationship between the max-weight update (\ref{eq:mw})-(\ref{eq:mw2}) and the Frank-Wolfe (FW) dual subgradient update (\ref{eq:fw})-(\ref{eq:fw2}).  The differences are that the optimisation (\ref{eq:mw}) is carried out over set $D$ (the set of \emph{actions} in max-weight terminology) rather than $C$ and approximate multiplier $\vv{\mu}_k$ is used instead of exact multiplier $\vv{\lambda}_{k}$.  Because $C:=\conv{(D)}$ then updates (\ref{eq:mw}) and (\ref{eq:fw}) are equivalent since the solution to a linear programme lies at an extreme point (and $D$ contains all the extreme points of $C$).  The substantive difference between the max-weight and FW updates therefore lies only in the use of $\vv{\mu}_k$ instead of $\vv{\lambda}_{k}$.

Regarding convergence analysis of the max-weight update (\ref{eq:mw})-(\ref{eq:mw2}), we can apply Theorem \ref{th:FWunsynch} with update set $\mathcal{U}$ containing the single element $u_1=\{1,\dots,n\}$ and using the Lagrangian $L(\cdot,\vv{\lambda}_k)$ as $F_k$.   Since the constraints are linear and the objective has bounded curvature the Lagrangian has bounded curvature (there is no need to restrict the magnitude of $\vv{\mu}_k$ when the constraints are linear).  Selecting step size $\alpha$ such that (\ref{eq:slow}) is satisfied and $\beta$ satisfying (\ref{eq:beta}) then, by Theorem \ref{th:FWunsynch}, we have that $\vv{z}_k\in C_\epsilon(\vv{\mu}_k)$ for $k$ sufficiently large and so condition (i) of  Theorem \ref{th:maintheorem} is satisfied.  By Lemma \ref{th:auxiliarymultiplier},  $\vv{\mu}_k$ satisfies condition (ii) of Theorem \ref{th:maintheorem} with $\sigma_0=2m(\sigma_1/\beta+\sigma_2)$, $\sigma_1=2\max_{\vv{z}\in C}\|\vv{A}\vv{z}\|_\infty$, and we are done. 

Observe that we can immediately extend the classical max-weight approach to allow nonlinear constraints by replacing update (\ref{eq:mw2}) with:
\begin{align}
&\vv{\mu}_{k+1} = [ \vv{\mu}_k  + \alpha \vv{g}_k(\vv{x}_k)]^{[0,\bar \lambda]}\label{eq:nl}
\end{align}
provided $\|\sum_{i=1}^k\vv{g}_i(\vv{x}_i) -\vv{g}(\vv{z}_i)\|_\infty \le \sigma_0/2$ for all $k=1,2,\dots$, where $\vv{g}_k:\mathbb{R}^n\rightarrow\mathbb{R}^m$, $k=1,2,\dots$.   
Note that we now need to use $\bar \lambda$ to place an upper bound on $\vv{\mu}_k$ in order to ensure that the Lagrangian has bounded curvature (as already noted, this is not needed when the constraints are linear).   By Lemma \ref{th:queuecontinuity}, $\|\sum_{i=1}^k\vv{g}_k(\vv{x}_i) -\vv{g}(\vv{z}_i)\|_\infty \le \sigma_0/2$ ensures that $\vv{\mu}_k$ satisfies condition (ii) of Theorem \ref{th:maintheorem}.   

Note that other nonlinear updates can also be used if more convenient than (\ref{eq:nl}).  For example, when $\vv{g}_k=\vv{g}$ using
\begin{align}\label{eq:gnl}
&\vv{\mu}_{k+1} = [ \vv{\mu}_k  + \alpha (\vv{g}(\vv{z}_k)+\partial\vv{g}(\vv{z}_k)^T (\vv{x}_k-\vv{z}_k))]^{[0,\bar \lambda]}
\end{align}
then
\begin{align*}
\|\vv{\mu}_{k+1}-\vv{\lambda}_k\|_\infty
& \le 2\alpha\max_{1\le j\le k} \Big\|\sum_{i=1}^j\partial\vv{g}(\vv{z}_i)^T (\vv{x}_i-\vv{z}_i) \Big\|_\infty \notag \\
& \le 2\alpha m \sigma_3\partial\bar{g}
\end{align*}
provided $\|\sum_{i=1}^k\vv{z}_i-\vv{x}_i\|_\infty\le \sigma_3$ $\forall k$ and $\partial\bar{g}:=\max_{\vv{z}\in C}\|\partial\vv{g}(\vv{z})\|_\infty$.  Since $\|\vv{\lambda}_{k+1}-\vv{\lambda}_k\|_\infty\le \alpha m\bar{g}$ then $\|\vv{\mu}_{k+1}-\vv{\lambda}_{k+1}\|_\infty \le \alpha\sigma_0$ as required by condition (ii) of Theorem \ref{th:maintheorem}, where $\sigma_0=m(\bar{g}+2 \sigma_3\partial\bar{g})$.


\subsubsection{Dual Max-Weight}
Consider the following update proposed in \cite{valls2014max} to solve convex optimisation $P^\prime$ with linear constraints:
\begin{align}
&\vv{x}_k\in\arg\min_{\vv{x}\in D} L((1-\beta)\vv{z}_k+\beta \vv{x},\vv{\mu}_k) \label{eq:pmw}\\
&\vv{z}_{k+1}=(1-\beta)\vv{z}_k+\beta \vv{x}_k\\
&\vv{\mu}_{k+1} = [ \vv{\mu}_k  + \alpha (\vv{A}\vv{x}_k-\vv{b}_k)]^+\label{eq:pmw2}
\end{align}
where $\|\sum_{i=1}^k\vv{b}_i -\vv{b}\|_\infty \le \sigma_2$ for all $k=1,2,\dots$.   Observe that this generalises  the classical dual subgradient update in a similar way that classical max-weight generalises the Frank-Wolfe dual subgradient update.
We can apply Theorem \ref{th:unsynch} with $\mathcal{U}$ containing the single element $u_1=\{1,\dots,n\}$ so that all elements of $\vv{z}_k$ are updated at each step and using the Lagrangian $L(\cdot,\vv{\lambda}_k)$ as $F_k$.   Provided step sizes $\alpha$, $\beta$ are selected sufficiently small then it follows that $\vv{z}_k\in C_\epsilon(\vv{\mu}_k)$ for $k$ sufficiently large and so satisfies Theorem \ref{th:maintheorem} condition (i).   
By Lemma \ref{th:auxiliarymultiplier},  $\vv{\mu}_k$ satisfies condition (ii) of Theorem \ref{th:maintheorem}. 

Similarly to classical max-weight, we can immediately extend the dual max-weight approach to include nonlinear constraints by replacing (\ref{eq:pmw2}) with update (\ref{eq:nl}) or (\ref{eq:gnl}).

\section{Generalising Max-Weight}

The analysis in the previous section makes clear that existing max-weight approaches use a specific running average method to ensure that approximate multiplier $\vv{\mu}_k$ remains close to the exact multiplier $\vv{\lambda}_k$ (and Lemma \ref{th:auxiliarymultiplier} establishes that a running average is indeed sufficient for this).  However, this naturally raises the question as to whether other approaches exist for choosing the sequences $\vv{z}_k$ and $\vv{x}_k$ such that  $\vv{\mu}_k$ and $\vv{\lambda}_k$ remain close.

One approach is to simply select $\vv{\mu}_k=\vv{\lambda}_k$ and then we are free to choose $\vv{z}_k$ and $\vv{x}_k$ as we like provided they ensure convergence to the $\epsilon$-ball $C_\epsilon(\vv{\mu}_k)$.   This is the approach taken by most (if not all) non-max-weight approaches that make use of multipliers.   However, it is clear from Lemma \ref{th:queuecontinuity} that any sequences $\vv{z}_k$ and $\vv{x}_k$ such that $\|\sum_{i=1}^k\vv{z}_i-\vv{x}_i\|_\infty$ is uniformly bounded for all $k=1,2,\dots$ will ensure that $\|\vv{\mu}_k-\vv{\lambda}_k\|_\infty$ is uniformly bounded\footnote{As already noted in Section \ref{sec:cont}, when $\vi{\lambda}{j}_k=0$ for all $k=1,2,\dots$ then it is sufficient that $\sum_{i=1}^k \vi{x}{j}_i$ is bounded above.  That is, for constraints which are not tight we have much greater flexibility.} 
and can be made as small as we like by adjusting the step size $\alpha$.   The following theorem shows that for any sequence $\vv{z}_k$, $k=1,2,\dots$ we can directly construct a sequence $\vv{x}_k$, $k=1,2,\dots$ that has this property.

\begin{theorem}[Discrete Actions]\label{th:discretesequence}
Let $D:=\{ \vv{x}_{1},\dots,\vv{x}_{|D|} \}$ be a subset of $\mathbb R^n$ and $C=\conv(D)$.  Let $\vv{z}_k$, $k=1,2,\dots$ be a sequence of points from $C$.  Since $C=\conv(D)$ then $\vv{z}_k=\vv{X}\vv{a}_k$ where matrix $\vv{X}=[\vv{x}_{1},\dots,\vv{x}_{|D|}]$ and $\vv{a}_k\in \Delta:=\{\vv{a}\in\mathbb{R}^{|D|}: \sum_{i=1}^{|D|}\vi{a}{i}=1, \vi{a}{i}\ge 0\}$.  Select $\vv{b}_{k+1} \in E$ such that
\begin{align}
\vv{S}_k+\vv{a}_{k+1}-\vv{b}_{k+1} \ge -\underline{S}\1\label{eq:bk}
\end{align}
where $\vv{S}_k=\sum_{i=1}^k\vv{a}_k-\vv{b}_k$, $E=\{\vv{e}_1,\dots,\vv{e}_{|D|}\}\subset \Delta$ and $\vv{e}_i$ denotes the unit vector with all elements zero apart from element $i$.   Provided $\underline{S}\ge 1$, then such a $\vv{b}_{k+1}$ always exists.  Further, the corresponding sequence $\vv{x}_k=\vv{X}\vv{b}_k$ of points from $D$ satisfies $\|\sum_{i=1}^k \vv{z}_i-\vv{x}_i\|_\infty \le (|D|-1)\underline{S}\|\vv{X}\|_\infty$ for all $k=1,2,\dots$. 
\end{theorem}
\begin{IEEEproof}
See the appendix.
\end{IEEEproof}

Note that many strategies exist to select $\vv{b}_{k+1} \in E$ satisfying (\ref{eq:bk}).  One is to select $\vv{b}_{k+1}$ according to
$
\vv{b}_{k+1}\in\arg\min_{\vv{b}\in E} \|\vv{S}_k+\vv{a}_{k+1}-\vv{b}\|_\infty
$.  
Another is to identify $\P_k:=\{j:\vi{S}{j}_k+\vi{a}{j}_{k+1}\ge 1-\underline{S}, j\in\{1,\cdots,n\}\}$, set $\P_k$ being non-empty since at least one such entry always exists by Theorem \ref{th:discretesequence}.  Selecting $\vv{b}_{k+1}=\vv{e}_j$ for some $j\in\P_k$ ensures $\vv{S}_k+\vv{a}_{k+1}-\vv{b}_{k+1} \ge -\underline{S}\1$.
%
%
%
Note also that when set $D$ is $\U$-feasible and so $C$ has product form $\prod_{u\in\U}C_u$ then selecting $\vv{x}_k$ can be decomposed into separate tasks of selecting each set of elements $\vvu{x}{u}_k$ such that $\|\sum_{i=1}^k \vvu{z}{u}_i-\vvu{x}{u}_i\|_\infty$ is bounded.

Theorem \ref{th:discretesequence} is constructive and provides a means for finding the $\vv{x}_k$, $k=1,2,\dots$.      It yields an approach where the $\vv{x}_k$'s are selected in an efficient online manner.  That is, given a new $\vv{z}_k\in C$ we can immediately determine a suitable $\vv{x}_k\in D$ via direct search.  The only state that needs to maintained in order to make this choice is a single vector $\vv{S}_k$ that keeps running track of the error.    Note that $\|\sum_{i=1}^k \vv{z}_i- \vv{x}_i \|_\infty  \le (|D|-1)\underline{S}\|\vv{X}\|_\infty$ implies that $\|\frac{1}{k}\sum_{i=1}^k \vv{z}_i- \vv{x}_i \|_\infty  \le \frac{1}{k}(|D|-1)\underline{S}\|\vv{X}\|_\infty$.  \dl{Hence, as $k\rightarrow\infty$ then $\frac{1}{k} \sum_{i=1}^k \vv{x}_i := \vv{x}^\diamond_k \rightarrow \vv{z}^\diamond_k$.  Since $f$ is convex (so Lipchitz continuous on compact set $C$) we also have that:
\begin{align}
|f(\vv{x}_k^\diamond) - f(\vv{z}_k^\diamond) | \le 
\frac{ (|D|-1)\underline{S}\|\vv{X}\|_\infty }{k} \nu_f
\end{align}
and so ${f}(\vv{x}_k^\diamond)\rightarrow {f}(\vv{z}_k^\diamond)$ as $k\rightarrow \infty$, where $\nu_f$ is the Lipschitz constant.} {Similarly, constraints $\vv{g}(\vv{x}_k^\diamond)\rightarrow \vv{g}(\vv{z}_k^\diamond)$ as $k\rightarrow \infty$.}

%

\begin{figure}
\centering
\subfloat[${z}_k$ and ${x}_k$]{\includegraphics[width=0.49\columnwidth]{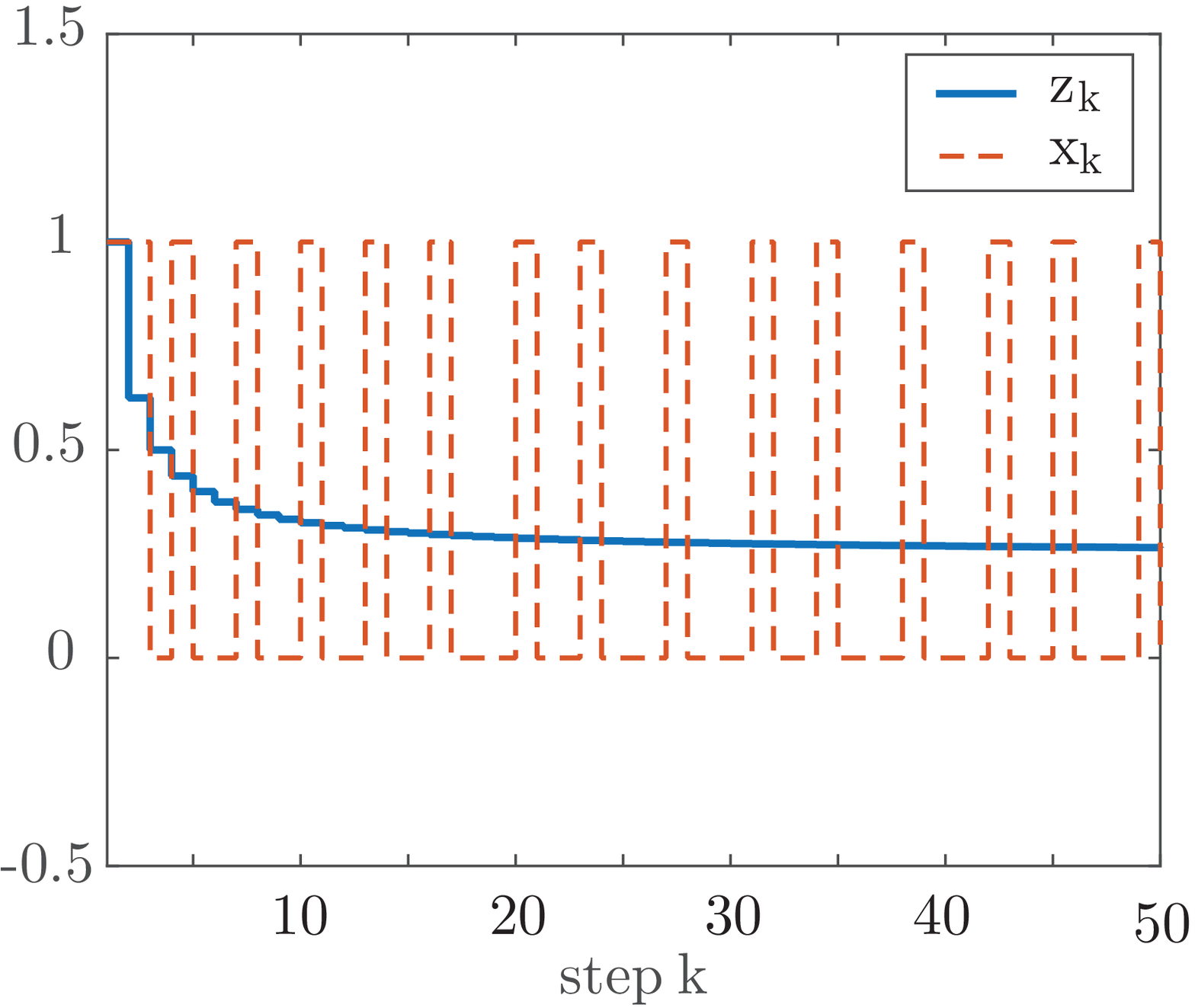}}
\hspace{1px}
\subfloat[$|\sum_{k=1}^K {z}_k- {w}_k |$]{\includegraphics[width=0.49\columnwidth]{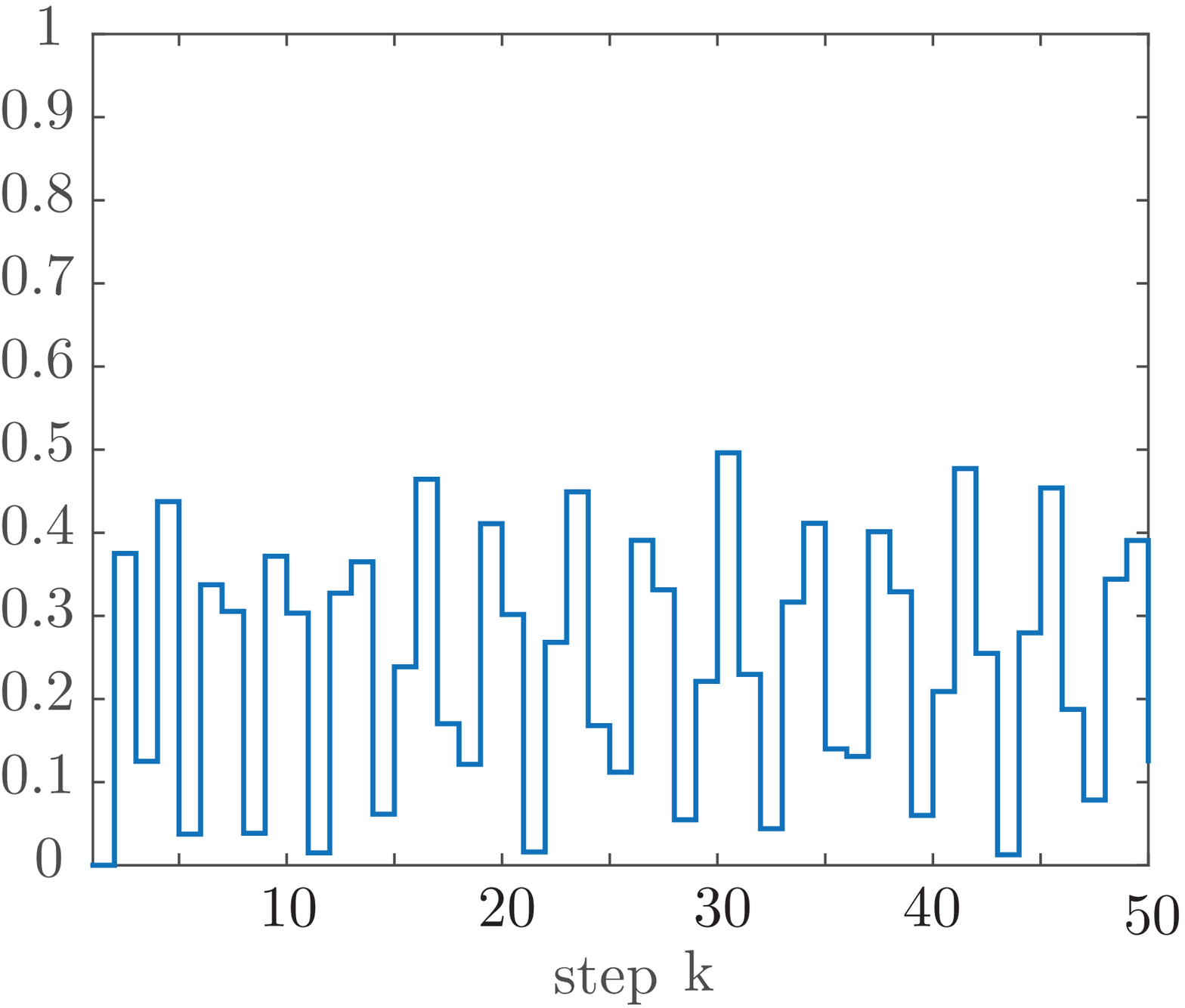}}
\caption{Illustrating sequence $x_k$ constructed according to Theorem \ref{th:discretesequence}.}\label{fig:ex2}
\end{figure}

Figure \ref{fig:ex2} illustrates construction of a sequence $x_k$ using the approach in Theorem \ref{th:discretesequence}.  In this example set $D=\{0,1\}$, $C=[0,1]$ and sequence $z_k = 0.75/k+0.25$, $k=1,2,\dots$.  It can be seen from Figure \ref{fig:ex2}(b) that $|\sum_{k=1}^K {z}_k- {w}_k |$ remains uniformly bounded for all $k=1,2,\dots$.     

Using Theorem \ref{th:discretesequence} for selecting $\vv{x}_k\in D$ given $\vv{z}_k\in C$ we can immediately extend essentially any existing multiplier-based approach to use discrete actions/queues.    

\subsection{Discrete Dual Approach}
For optimisation $P^\prime$, using Theorem \ref{th:discretesequence} the classical dual subgradient approach can be directly generalised to:
\begin{align}
&\vv{z}_{k} \in \arg\min_{\vv{z}\in C} L(\vv{z},\vv{\mu}_k) \label{eq:pm4}\\
&\text{Select$\footnotemark$  $ \vv{x}_{k}\in D$ s.t. $\Big\|\sum_{i=1}^k \vv{z}_i - \vv{x}_i \Big\|_\infty 
 \le \sigma_3$}\\
&\vv{\mu}_{k+1} = [ \vv{\mu}_k  + \alpha (\vv{A}\vv{x}_k-\vv{b}_k)]^+\label{eq:pm5}
\end{align}
\footnotetext{When $\vi{\lambda}{j}_k=0 \ \forall k$ then it is sufficient that $\sum_{i=1}^k \vi{x}{j}_i$ is bounded above.}
where $\|\sum_{i=1}^k\vv{b}_i -\vv{b}\|_\infty \le \sigma_2$ for all $k=1,2,\dots$.  As already noted, this choice of $\vv{z}_k$ satisfies Theorem \ref{th:maintheorem} condition (i), namely $\vv{z}_k \in  C_\epsilon (\vv{\mu}_k)$ with $\epsilon=0$.  It follows from Lemma \ref{th:queuecontinuity} that when $\|\sum_{i=1}^k \vv{z}_i - \vv{x}_i \|_\infty 
 \le \sigma_3$ then $\| \vv{\lambda}_k - \vv{\mu}_k \|_\infty \le 2\alpha \sigma_0$ with $\sigma_0=m(\max_{\vv{z}\in C}\|\vv{A}\vv{z}\|_\infty\sigma_3+\sigma_2)$, and so Theorem \ref{th:maintheorem} condition (ii) is also satisfied.

Similarly to before, update (\ref{eq:pm4})-(\ref{eq:pm5}) can be further extended to allow use of unsynchronised updates provided the Lagrangian is $\U$-separable, action set $D$ is $\U$-feasible and $\U$ is an admissible update set.  Namely,
\begin{align}
&\vvu{z}{u}_{k} \in \arg\min_{\vv{w}\in C_u} \vi{L}{u}(\vv{w},\vvi{\mu}{u}_k),\ u\in\U \label{eq:upm4}\\
&\text{Select$^5$  $\vv{x}_{k}\in D$ s.t. $\Big\|\sum_{i=1}^k \vv{z}_i - \vv{x}_i \Big\|_\infty   \le \sigma_3$}\\
&\vv{\nu}_{k+1} = [ \vv{\nu}_k  + \alpha (\vv{A}\vv{x}_k-\vv{b}_k)]^+\label{eq:upm5}
\end{align}
where $C = \prod_{u\in\U}C_u$ with $C_u\subset \mathbb{R}^{|u|}$. The elements $u$ which are updated use multiplier $\vvi{\mu}{u}_k=\vv{\nu}_{k}$.  Those elements $\vvu{z}{v}_k$, $v\ne u$ which are not updated at step $k$ use the old multiplier value $\vvi{\mu}{v}_k=\vv{\nu}_{k-\vi{\tau}{v}_k}$ where $k-\vi{\tau}{v}_k$ is the time step where $\vvu{z}{v}_k$ was last updated. 

Consideration can be extended to nonlinear constraints by replacing update (\ref{eq:upm5}) with (\ref{eq:gnl}) or (\ref{eq:nl}) (in the latter case selecting the $\vv{x}_k$, $k=1,2,\dots$ such  $\|\sum_{i=1}^k\vv{g}_i(\vv{x}_i) -\vv{g}(\vv{z}_i)\|_\infty \le \sigma_0/2$ for all $k=1,2,\dots$).

Note that the approach here does not require the objective or constraints to have bounded curvature.   This is in marked contrast to most existing max-weight approaches, and so updates (\ref{eq:pm4})-(\ref{eq:pm5}) and (\ref{eq:upm4})-(\ref{eq:upm5}) can be applied to a much wider class of problems.  

Observe also that the unconstrained convex optimisation in update (\ref{eq:upm4}) can be solved in many ways (gradient-descent, proximal, Newton-based, etc.) with the most appropriate method depending on the properties of the Lagrangian.  Further, by Theorem \ref{th:maintheorem}, convergence is only needed to an $\epsilon$-ball $C_\epsilon(\vv{\mu}_k)$ about the optimum.


\subsection{Example: Link With Communication Delays} \label{sec:link_example}

\begin{figure}
\centering
\includegraphics[width=0.6\columnwidth]{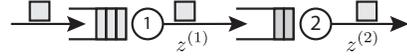}
\caption{Schematic illustrating example in Section \ref{sec:link_example}.}
\label{fig:link_example}
\end{figure}

Consider the link illustrated in Figure \ref{fig:link_example}. Packets arrive in the queue of node $1$, node $1$ transmits packets to node $2$, and the packets transmitted by node $2$ leave the system. Time is slotted and for simplicity we will assume that slots have fixed duration and that nodes can transmit only one packet in a slot. 

Our aim is to make scheduling decisions at each node (namely, transmit or not transmit) that minimise a convex  function of the average throughput.  Action set $D=\{0,1\}^2$, where a $1$ in element $\vi{x}{i}$ of $\vv{x}\in D$ indicates a packet transmission from queue $i$.  Update set $\U=\{u,v\}$ with $u=\{1\}$, $v=\{2\}$.  The convex formulation of the problem is: 
\begin{align}
\min_{\vv{z} \in C} \qquad & f(\vv{z}):=\vi{f}{u}(\vi{z}{1}) + \vi{f}{v}(\vi{z}{2}) \\
\text{s.t.}  \qquad  & b - \vi{z}{1} \le 0,\  \vi{z}{1} - \vi{z}{2} \le 0
\end{align}
where $C =  \conv{(D)}$ and $b$ the average number of packets that arrive into the queue of node $1$ in a slot. The Lagrangian is given by $L(\vv{z}, \vv{\lambda}) = \vi{L}{u}(\vi{z}{1},\vv{\lambda}) + \vi{L}{v}(\vi{z}{2},\vv{\lambda})$ where $\vi{L}{u}(z,\vv{\lambda}) = \vi{f}{u}(z) + (\vi{\lambda}{2} - \vi{\lambda}{1}) z+ \vi{\lambda}{1}b$, $\vi{L}{v}(z,\vv{\lambda}) = \vi{f}{v}(z) - \vi{\lambda}{2} z$ with $\vv{\lambda} = [\vi{\lambda}{1},\vi{\lambda}{2}]^T$, {i.e.}, the Lagrangian is $\U$-separable.  We can use the discrete dual approach in (\ref{eq:upm4})-(\ref{eq:upm5}) to solve this problem in a distributed manner. In particular, at each time slot each node solves the convex optimisation
\begin{align}
\dynwidth{
\vi{z}{1}_k \in \underset{w \in [0,1]}{\arg \min} \ \vi{L}{u}(w, \vvi{\mu}{u}_k),
\vi{z}{2}_k \in \underset{w \in [0,1]}{\arg \min} \ \vi{L}{v}(w, \vvi{\mu}{v}_k) 
}\label{eq:2zupdate}
\end{align}
and chooses discrete action,
\begin{align}
\vi{x}{i}_k  \in \underset{x \in \{0,1\}}{\arg \min} \ \| (\textstyle \sum_{j=1}^k \vi{z}{i}_j - \textstyle \sum_{j=1}^{k-1} \vi{x}{i}_j) - x \|_\infty  \label{eq:2xupdate}
\end{align}
where $\vvi{\mu}{u}_k = \alpha [ \vi{Q}{1}_k , \vi{Q}{2}_{k - \vi{\tau}{1}}]^T$, $\vvi{\mu}{v}_k = \alpha [ \vi{Q}{1}_{k-\vi{\tau}{2}} , \vi{Q}{2}_{k}]^T$ with queue updates:
\begin{align*}
\vi{Q}{1}_{k+1} & = [ \vi{Q}{1}_k - \vi{x}{1}_k + b_k]^+,\ \vi{Q}{2}_{k+1} = [ \vi{Q}{2}_k - \vi{x}{2}_k + \vi{x}{1}_k]^+.
\end{align*}
Observe that the elements of $\vvi{\mu}{u}_k$, $\vvi{\mu}{v}_k$ capture delay in the exchange of network state information via parameters $\vi{\tau}{1}$ and $\vi{\tau}{2}$.   We assume that $\vi{\tau}{i} \in \{0, 1,\dots, \bar \tau\}, \ i=1,2$ and $\bar \tau$ is a constant.  Despite being delayed, by Theorem \ref{th:discretesequence} the scaled queue occupancies $\vvi{\mu}{u}_k$, $\vvi{\mu}{v}_k$ still approximate the exact multipliers $\vv{\lambda}_k$. 
Our choice of approximate multipliers in this example only captures delays in the exchange of network state information, however, by selecting them appropriately we could also consider asynchronous updates in (\ref{eq:2zupdate}).

\begin{figure}
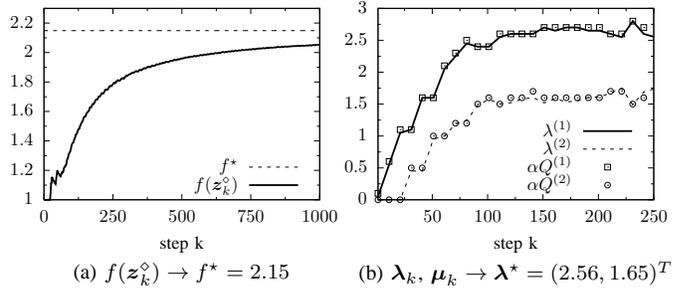

\centering
\subfloat[$f(\vv{z}^\diamond_k) \rightarrow f^\star = 2.15$]
{
\resizebox{0.46\columnwidth}{!}{\input{link_f.tex}}
}\label{fig:link_f}
\subfloat[$\vv{\lambda}_k$, $\vv{\mu}_k \to \vv{\lambda}^\star = (2.56,1.65)^T$ ]
{
\resizebox{0.46\columnwidth}{!}{\input{link_q.tex}}
}\label{fig:link_q}
\caption{Illustrating discrete dual subgradient solution to link example (Section \ref{sec:link_example})}\label{fig:link}
\end{figure}

We simulate the network with utility functions\footnote{$\vi{f}{v}$ does not have bounded curvature.} $\vi{f}{u}(z) = z$, $\vi{f}{v}(z) = \max \{ \exp(z), \pi z\}$, step size $\alpha = 0.1$ and $\{ b_k, k=1,2,\dots, K \}$ a $\{0,1\}$ sequence that satisfies $|\sum_{i=1}^k (b_k - b) | \le \frac{1}{2}$ where $b = \frac{1}{2}$. Also, we use initial condition $\vv{\lambda}_1 = \vv{Q}_1 = \vv{0}$ and let $\vi{\tau}{i}, i=1,2$ be a random variable that takes values uniformly at random between $\{ 0, \dots,  \bar \tau\}$, $\bar \tau = 5$. Figure \ref{fig:link}(a) shows the convergence of $f(\vv{z}^\diamond_k)$ to $f^\star$, and Figure \ref{fig:link}(b) shows how the Lagrange multipliers and the $\alpha$-scaled queues converge to a ball around $\vv{\lambda}^\star  = [ 2.56 ,  1.65]^T$. Importantly, we never require to know the exact Lagrange multipliers in the optimisation and in Figure \ref{fig:link}(b) we compute them for illustrative purposes only. 

Important features illustrated by this example include (i) use of scaled queues occupancies as surrogates for the Lagrange multipliers in the optimisation, (ii) nodes make decisions in a distributed manner and (iii) the dual subgradient update makes use of imperfect network state information yet still converges. The distributed aspect ensures that the complexity of finding an optimal scheduling does not grow with the number of nodes in the network. 

%

We can easily extend this example to allow multiple packets to be transmitted in a single slot by defining ${{D = \{ 0, 1,\dots, s \}^2}}$ where $s$ denotes the maximum number of packets that can be transmitted in a slot. In effect, we are dividing a slot in $s$ \emph{subslots} but there is no constraint on the order in which packets have to be transmitted. For example, if $s=5$ and $x_k = 3$ transmitting the sequence of packets $\{ 1,1,1,0,0\}$ in a slot is equivalent to transmitting sequence $\{1,0,1,0,1\}$ since all it matters is that $3$ packets are transmitted. Note that queue continuity trivially holds. With this setup in each time slot we are deciding how many (not how) packets to transmit in the $s$ following subslots, i.e.,  slots can be regarded as \emph{refresh times}.
Since the duration of a slot is fixed, by selecting $s$ large we can change packets for bits in the problem setup, and by adding constraints/requirements on how a stream of bits can be transmitted we are indirectly considering the case where packet transmissions have different durations.  

\subsection{Generalising Classical Max-Weight}

We can generalise the Frank-Wolfe subgradient approach to allow use of unsynchronised updates {i.e.}, only a subset of elements are updated at each time step.  By Theorem \ref{th:FWunsynch} we can still ensure convergence and satisfaction of condition (i) of Theorem \ref{th:maintheorem} provided step sizes $\alpha$ and $\beta$ are sufficiently small.  

Unfortunately, when updates are unsynchronised then we can no longer use Lemma \ref{th:auxiliarymultiplier} to ensure that condition (ii) of Theorem \ref{th:maintheorem} is satisfied.   That is, a running average of the form $\vv{z}_{k+1}=\vv{z}_k+\beta \vv{U}_{u}(\vv{x}_k-\vv{z}_k)$ is not enough to ensure that approximate multiplier $\vv{\mu}_k$ stays close to the exact multiplier $\vv{\lambda}_k$.   However, we can use Theorem \ref{th:discretesequence} to ensure this.  That is, by using the following update for solving optimisation $P^\dagger$ with linear constraints:
\begin{align}
&\vv{y}_k \in\arg\min_{\vv{x}\in D} \partial_{\vv{z}}L(\vv{z}_k ,\vv{\mu}_k)^T \vv{U}_{u_{k}}\vv{x} \label{eq:umw}\\
&\vv{z}_{k+1}=\vv{z}_k+\beta \vv{U}_{u_{k}}(\vv{y}_k-\vv{z}_k)\label{eq:umw1}\\
&\text{Select$^5$ $\vv{x}_{k}\in D$ s.t. $\|\sum_{i=1}^k \vv{z}_i - \vv{x}_i \|_\infty 
 \le \sigma_3$}\label{eq:umw1b}\\
&\vv{\mu}_{k+1} = [ \vv{\mu}_k  + \alpha (\vv{A}\vv{x}_k-\vv{b}_k)]^+
 \label{eq:umw2}
\end{align} 
where $\|\sum_{i=1}^k\vv{b}_i -\vv{b}\|_\infty \le \sigma_2$ for all $k=1,2,\cdots$ and $u_k\in\U$, $k=1,2,\dots$ is a sequence that regularly visits every element of $\vv{z}$.   

Selecting step size $\alpha$ such that (\ref{eq:slow}) is satisfied and $\beta$ satisfying (\ref{eq:beta}) then applying Theorem \ref{th:FWunsynch} using the Lagrangian $L(\cdot,\vv{\lambda}_k)$ as $F_k$ we have that $\vv{z}_k\in C_\epsilon(\vv{\mu}_k)$ for $k$ sufficiently large and so condition (i) of  Theorem \ref{th:maintheorem} is met.   By Lemma \ref{th:queuecontinuity}, $\vv{\mu}_k$ satisfies condition (ii) of Theorem \ref{th:maintheorem}.

As usual consideration can be extended to nonlinear constraints by replacing update (\ref{eq:umw2}) with (\ref{eq:gnl}) or, alternatively, by using update (\ref{eq:nl}) and selecting the $\vv{x}_k$, $k=1,2,\dots$ such  $\|\sum_{i=1}^k\vv{g}_i(\vv{x}_i) -\vv{g}(\vv{z}_i)\|_\infty \le \sigma_0/2$ for all $k=1,2,\dots$.

Recall that coordinate descent is recovered as a special case when {e.g.}, $\U=\{\{1\},\dots,\{n\}\}$ and $D=\{\vv{e}_1,\dots,\vv{e}_n\}\cup \{\1\}$.  
However, update (\ref{eq:umw})-(\ref{eq:umw2}) is a significant generalisation beyond this.


\subsection{Example: Unsynchronised Max-Weight Updates}\label{sec:ex6}

Consider a collection of $n$ queues.   Jobs arrive at queue $i$ at mean rate $\vi{b}{i}$ and are buffered until they can be served.   Up to $d$ jobs can be served from a queue in a time slot, but this incurs a cost proportional to the square of the number of jobs processed at one time.     In action vector $\vv{x}$ element $\vi{x}{i}=1$ indicates service of one job from queue $i$, $\vi{x}{i}=2$ service of 2 jobs and so on.   The set of actions is $D=\{0,\dots,d\}^n$.   The objective is to minimise $f(\vv{z}) = \sum_{i=1}^n(\vi{z}{i})^2$ subject to the constraints that $\vi{b}{i}\le \vi{z}{i}$ for $i=1,\dots,n$ (the service rate at queue $i$ is no less than the arrival rate).   

At each time slot an element of $\U=\{\{1\},\{2\},\dots,\{n\}\}$ is selected cyclically and the corresponding element of vector $\vv{z}_k$ updated according to (\ref{eq:umw})-(\ref{eq:umw1}).  Note that this means each update is simple and fast as it only involves deciding which of $d+1$ values to use, independent of $n$.

\begin{figure}
\centering
\subfloat[$f(\vv{z}^\diamond_k)\rightarrow f^*$]{\includegraphics[width=0.46\columnwidth]{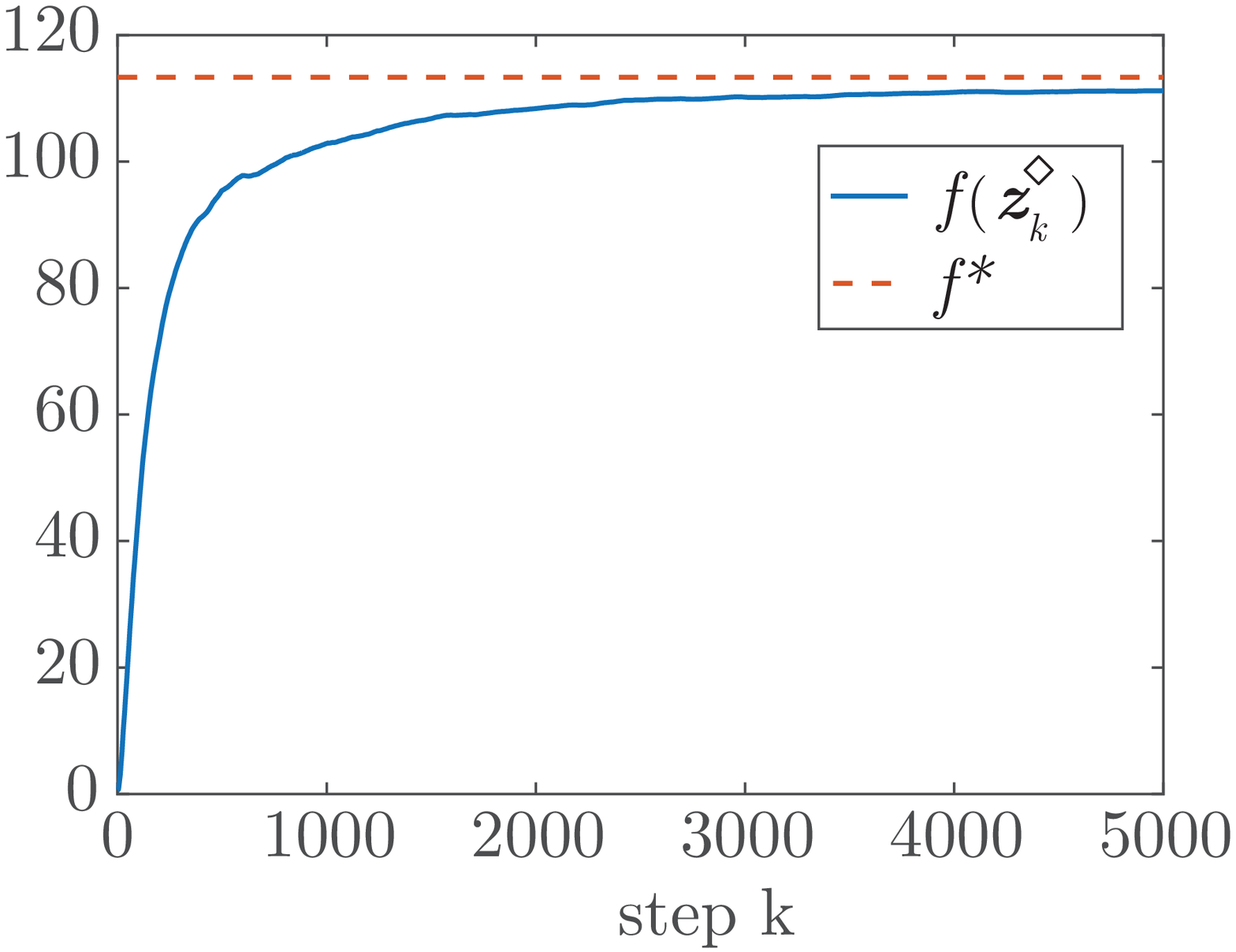}}
\hspace{1px}
\subfloat[$\vv{\mu}_k \to \vv{\lambda}^\star$]{\includegraphics[width=0.46\columnwidth]{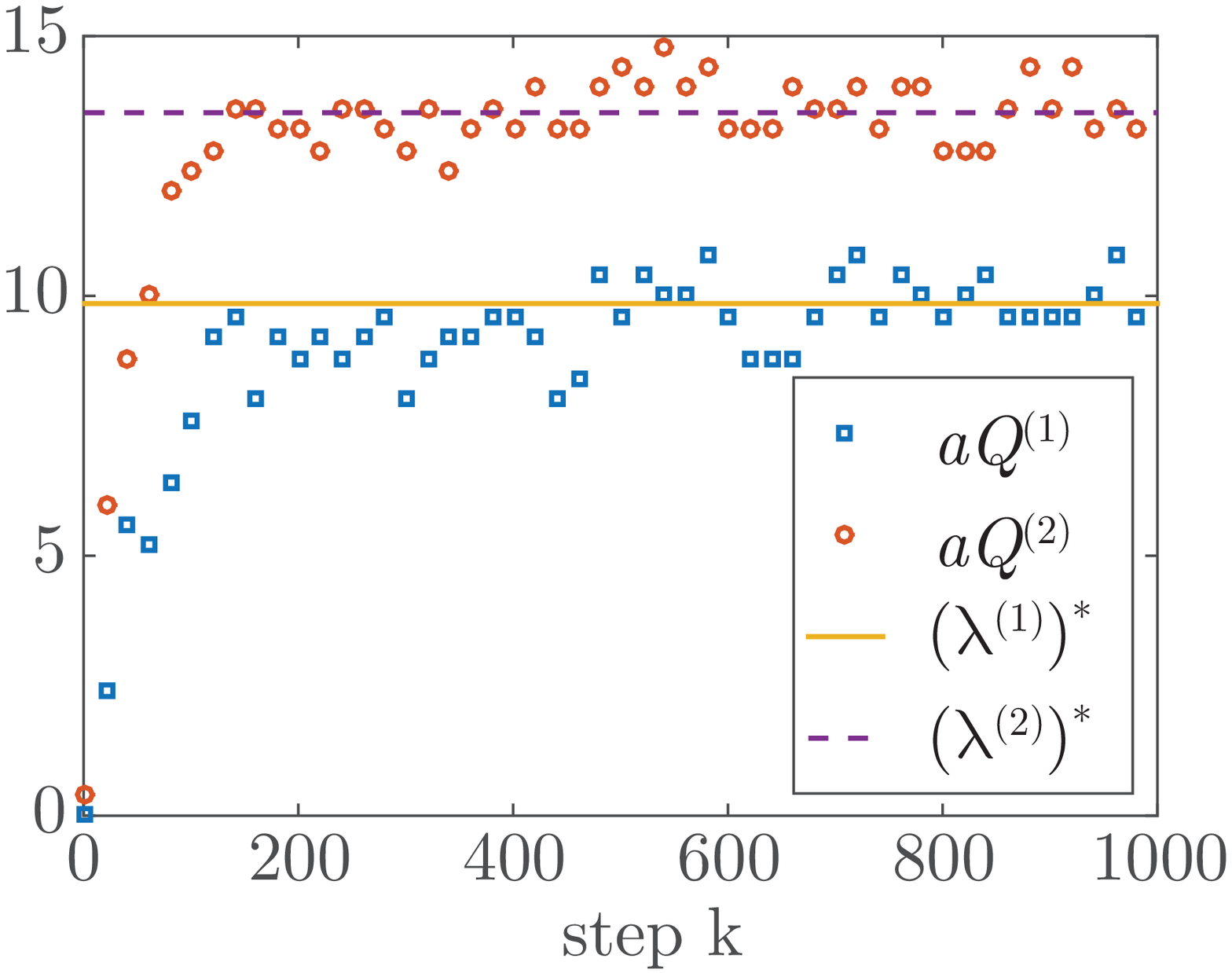}}
\caption{Illustrating application of unsynchronised max-weight update (\ref{eq:umw})-(\ref{eq:umw2}) to the example in Section\ref{sec:ex6}.}\label{fig:ex6}
\end{figure}

Figure \ref{fig:ex6}(a) illustrates the convergence of the objective to a ball around the optimum when $n=2$, $d=8$, step sizes $\alpha=0.05$ and $\beta=0.1$.  Figure \ref{fig:ex6}(b) shows how the corresponding time histories of the scaled queue occupancy $\vv{\mu}_k$ converge to a ball around the optimum multipliers $\vv{\lambda}^\star$.   

\subsection{Generalising Dual Max-Weight}
Similarly, we can extend the Dual Max-Weight approach to use unsynchronised updates for solving optimisation $P^\dagger$:
\begin{align}
&\vv{y}_k \in\arg\min_{\vv{x}\in D} L(\vv{z}_k+ \beta\vv{U}_{u_k} (\vv{x}-\vv{z_k}),\vv{\mu}_k)\label{eq:udm}\\
&\vv{z}_{k+1}=\vv{z}_k+ \beta \vv{U}_{k} (\vv{y}_k-\vv{z}_k)\\
&\text{Select$^5$ $\vv{x}_{k}\in D$ s.t. $\Big\|\sum_{i=1}^k \vv{z}_i - \vv{x}_i \Big\|_\infty 
 \le \sigma_3$}\\
&\vv{\mu}_{k+1} = [ \vv{\mu}_k  + \alpha(\vv{A}\vv{x}_k-\vv{b}_k)]^+
\end{align}
where $\|\sum_{i=1}^k\vv{b}_i -\vv{b}\|_\infty \le \sigma_2$ $\forall k$ and $u_k\in\U$, $k=1,2,\dots$ is a sequence that regularly visits every element of $\vv{z}$.   Applying Theorem \ref{th:unsynch} using the Lagrangian $L(\cdot,\vv{\lambda}_k)$ as $F_k$, provided $\alpha$ and $\beta$ are sufficiently small then $\vv{z}_k\in C_\epsilon(\vv{\mu}_k)$ for $k$ sufficiently large and so condition (i) of  Theorem \ref{th:maintheorem} is satisfied.  By Lemma \ref{th:queuecontinuity},  $\vv{\mu}_k$ satisfies condition (ii) of Theorem \ref{th:maintheorem}.

Once again, nonlinear constraints can be included by replacing update (\ref{eq:umw2}) with (\ref{eq:gnl}) or (\ref{eq:nl}), and coordinate descent corresponds to selecting $\U=\{\{1\},\dots,\{n\}\}$.   
Note also that while action set $D$ need consist only of the extreme points of set $C$, we can optionally select $D$ equal to $C$.   The optimisation (\ref{eq:udm}) is then convex and encompasses a line search for the optimum step size.

\subsection{Decisions \& Actions At Different Time Scales}

\begin{figure}
\centering
\subfloat[]{\includegraphics[width=0.49\columnwidth]{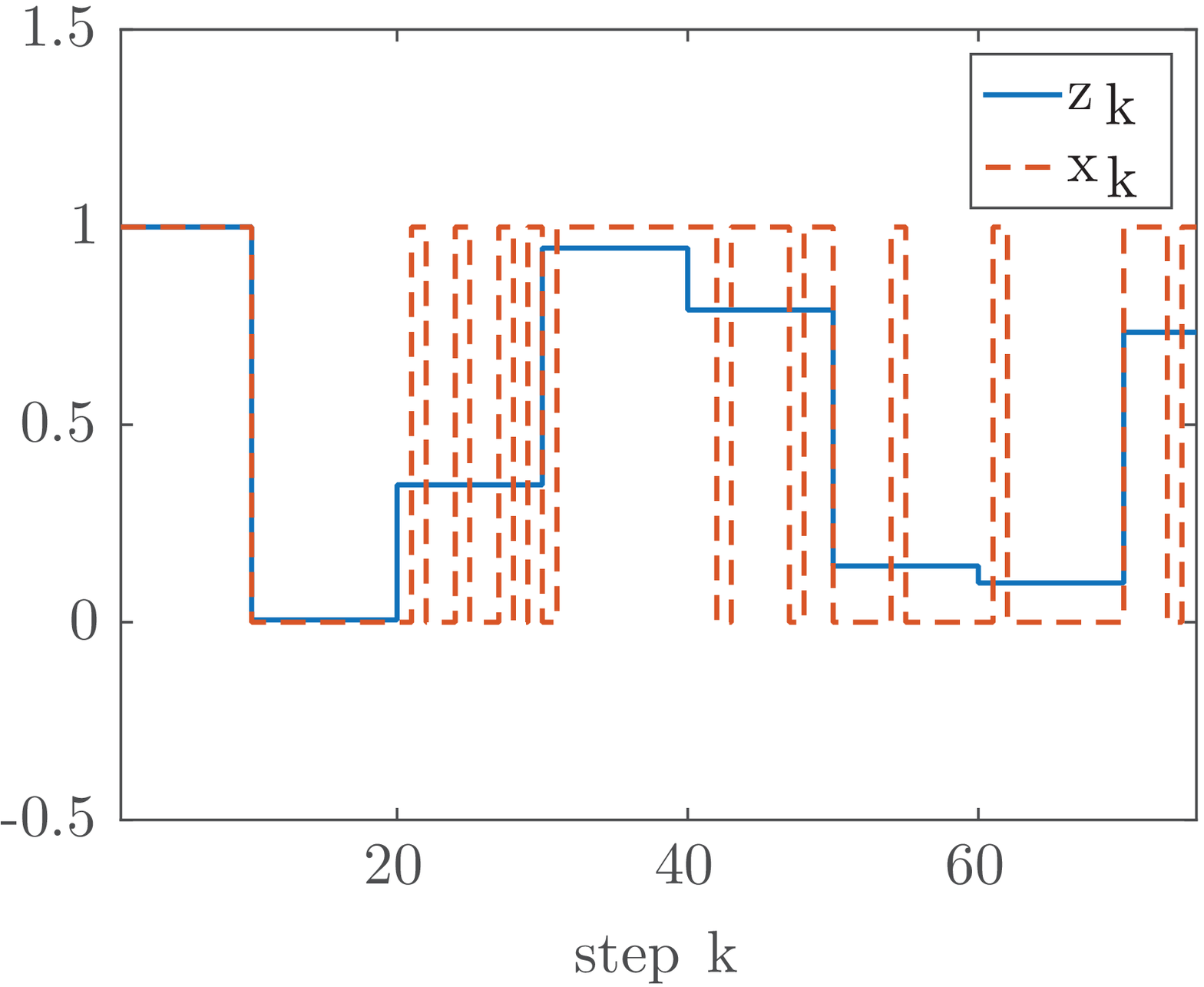}}
\hspace{1px}
\subfloat[]{\includegraphics[width=0.49\columnwidth]{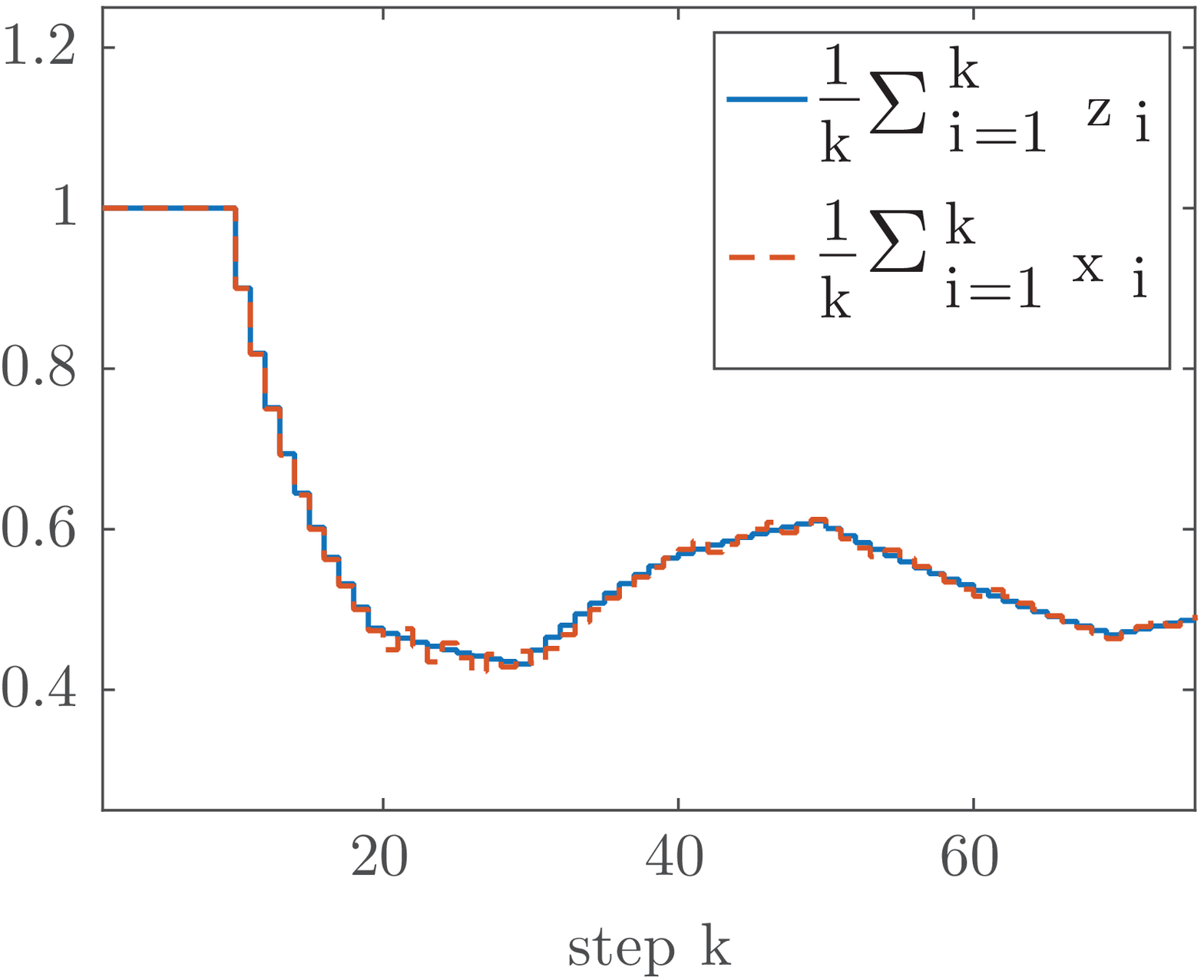}}
\caption{Illustrating sequence of actions $x_k$ updated on a faster time scale than the sequence of decisions $\vv{z}_k$. Although updated at different time scales, by Theorem \ref{th:discretesequence} we can still ensure that $\sum_{i=1}^kz_i$ and $\sum_{i=1}^kx_i$ stay close, see plot (b).}\label{fig:ex4}
\end{figure}

By decoupling the choice of $\vv{z}_k$ from the choice of $\vv{x}_k$ (in contrast to when the running average approach is used), Theorem \ref{th:discretesequence} also allows actions to be taken at a faster rate than that at which decisions are calculated.   This is illustrated for example in Figure \ref{fig:ex4} where $\vv{z}_k$ is updated every $10$ slots while $\vv{x}_k$ is updated every slot.  Since Theorem \ref{th:discretesequence} holds for all sequences $\vv{z}_k$, we can still select actions $\vv{x}_k$ that ensure $\|\sum_{i=1}^k\vv{z}_i- \vv{x}_i\|_\infty$ stays bounded.   This is of considerable practical importance in situations where actions $\vv{x}_k$ must be selected quickly, {e.g.}, selecting whether to transmit a packet or not needs to be taken every $8 \ \mu$s on a 1 Gbps link with 1000B packets, yet decisions  $\vv{z}_k$ take a longer time to be calculated, {e.g.}, updating $\vv{z}_k$ takes $10 \ ms$. 

\section{Delayed Communication \& Unsynchronised Multiplier Updates}
\subsection{Delayed Communication}
Condition (ii) of Theorem \ref{th:maintheorem} requires $\| \vv{\lambda}_k - \vv{\mu}_k \|_\infty \le \alpha \sigma_0$.   The slack between $\vv{\lambda}_k$ and $\vv{\mu}_k$ can accommodate communication delays and use of stale multiplier values.  
For example, suppose that the $i$'th element of approximate multiplier $\vv{\mu}$ has value equal to the $i$'th element of the exact multiplier $\vv{\lambda}$ from $\vi{\tau}{i}$ time steps in the past, that is $\vvi{\mu}{i}_k = \vvi{\lambda}{i}_{k-\vi{\tau}{i}}$.   Suppose also the delay is upper bounded, $\vi{\tau}{i}\le \bar{\tau}$, $i=1,\dots,m$.   Since $\vv{\lambda}_{k+1} = [ \vv{\lambda}_k  + \alpha \vv{g}(\vv{z}_k)]^+$ we have that $|\vi{\lambda}{i}_{k} - \vi{\lambda}{i}_{k-\tau}|\le \alpha \bar{g}\bar{\tau}$ and so $\| \vv{\lambda}_k - \vv{\mu}_k \|_\infty \le \alpha \sigma_0$, $k=1,2,\dots$  with $\sigma_0=\bar{g}\bar{\tau}$.  That is, condition (ii) of Theorem \ref{th:maintheorem} is still satisfied.


\subsection{Unsynchronised Multiplier Updates}
The slack between $\vv{\lambda}_k$ and $\vv{\mu}_k$ can also more generally accommodate unsynchronised updates to $\vv{\mu}_k$.
For example, suppose a packet is dequeued and transmitted but that due to propagation delay arrives at its destination queue some time later (this may be, for example, around 100 $ms$ later if the link is across the Atlantic).  Suppose the departure queue corresponds to element $i$ of $\vv{\mu}$ and the arrival queue corresponds to element $j$ of $\vv{\mu}$.   Then $\vi{\mu}{i}_k$ is updated when transmission starts and after delay $\tau$ then $\vi{\mu}{j}_{k+\tau}$ is updated.   These updates correspond to the selection of a single action, namely transmission of a packet across the link, yet in this case the effect of this action on $\vv{\mu}$ is not confined to the time slot in which the action is selected.   Nevertheless, when the update to $\vi{\mu}{j}_{k+\tau}$ happens delay $\tau$ after the exact multiplier $\vv{\lambda}_k$ is updated then, provided $\tau$ is bounded, by the same argument as before $\| \vv{\lambda}_k - \vv{\mu}_k \|_\infty$ can be kept small for all $k=1,2,\cdots$ and condition (ii) of Theorem \ref{th:maintheorem} satisfied.


\section{Stochasticity}
The analysis above is for deterministic optimisation problems.   However, it can be readily extended to a class of stochastic optimisations.   

\subsection{Approximate Multipliers}
\subsubsection{Queue Arrivals} Suppose we have linear constraints $\vv{A}\vv{z}-\vv{b}\preceq \vv{0}$ and use $\vv{A}\vv{x}_k-\vv{b}_k$ in the update of approximate multiplier $\vv{\mu}_k$.   This might correspond to a queueing network where $\vv{A}$ defines the interconnections between queues and $\vv{b}_k$ the arrivals.
Let $\{\vv{B}_k\}$ be a stochastic process with realisations taking values in $\mathbb{R}^m$ and mean $\vv{b}$. Let $p_K:=\Prob(\max_{ k\in\{1,2,\dots,K\}}\|\sum_{i=1}^k (\vv{B}_i -\vv{b})\|_{\infty} \le {\sigma_2})$. Let  $\{\vv{b}_i\}_{i=1}^K$ denote a  realisation of length $K$. Fraction $p_K$ of these realisations satisfy $\|\sum_{i=1}^k (\vv{b}_i-\vv{b})\|_{\infty} \le {\sigma_2}$ for all $k=1,2,\dots,K$.  When this fraction is asymptotically lower bounded $\liminf_{K\rightarrow\infty}p_K\ge p$, then our analysis applies to fraction $p$ of realisations {i.e.}, with probability $p$.
Note that there is no requirement for stochastic process $\{\vv{B}_k\}$ to be i.i.d. or for any of its properties, other than that feasible set $\vv{A}\vv{z}\preceq\vv{b}$ has non-empty relative interior, to be known in advance. 
Extension to nonlinear constraints follows similarly.  

\subsubsection{Actions} Our analysis applies with probability $p$ when sequences $\vv{x}_k$ are selected so as to ensure $\|\sum_{i=1}^k \vv{z}_i - \vv{x}_i \|_\infty$ is bounded with probability $p$.   This means that the $\vv{x}_k$ may, for example, be constructed using stochastic approaches such as Gibbs sampling.

\subsection{Actions With Random Effect}
Suppose that when at time $k$ we select action $\vv{x}_k\in D$, the action actually applied is a realisation of random variable $\vv{Y}_k$ that also takes values in $D$; this is for simplicity, the extension to random action sets different from $D$ is straightforward.  For example, we may select $x_k=1$ (which might correspond to transmitting a packet) but with some probability actually apply $y_k=0$ (which might correspond to a transmission failure/packet loss).  

Let $p_{\vv{x}\vv{y}}:=\Prob(\vv{Y}_k=\vv{y} | \vv{x}_k=\vv{x})$, $\vv{x}$, $\vv{y}\in D$ and we assume that this probability distribution is time-invariant \emph{i.e.,} does not depend on $k$; again, this can be relaxed in the obvious manner.    Then $\bar{\vv{y}}(\vv{x}) := \mathbb{E}[\vv{Y}_k|\vv{x}_k=\vv{x}] = \sum_{\vv{y}\in D} \vv{y} p_{\vv{x}\vv{y}}$ can be calculated.   The above analysis now carries over unchanged provided we everywhere replace $\vv{x}_k$ by $\bar{\vv{y}}(\vv{x}_k)$ {e.g.}, modify the non-convex optimisation from $\min_{\vv{x} \in D} L((1-\beta)\vv{z}_k + \beta \vv{x},{\vv{\lambda}}_k)$ to $ \min_{\vv{x} \in D} L((1-\beta)\vv{z}_k + \beta \bar{\vv{y}}(\vv{x}),{\vv{\lambda}}_k)$.  That is, we simply change variables to $\bar{\vv{y}}$.  

\section{Conclusions}
We study the use of approximate Lagrange multipliers and discrete actions in solving convex optimisation problems.   We show that descent, which can  be ensured using a wide range of approaches (subgradient, proximal, Newton, \emph{etc}.), is orthogonal to the choice of multipliers.   Using the Skorokhod representation for a queueing process we show that approximate multipliers can be constructed in a number of ways.  This leads to the generalisation of (i) essentially any descent method to encompass use of discrete actions and queues and (ii) max-weight scheduling to encompass new descent methods including those with unsynchronised updates such as block coordinate descent. 

\section*{Appendix: Proofs}

\subsection{ Descent}

\subsection*{Proof of Lemma \ref{th:coord}}

%
By convexity,
$F(\vv{z}) +  \partial F(\vv{z})^T\vv{U}_{u}(\vv{y}-\vv{z}) \le F(\vv{z}+\vv{U}_{u}(\vv{y}-\vv{z}))
 \le F(\vv{z})-\epsilon$.
Hence, $\partial F(\vv{z})^T\vv{U}_{u}(\vv{y}-\vv{z}) \le -\epsilon$.    Since $D$ is $\U$-feasible, for any $\vv{x}\in D$ we have $ \vv{z} + \beta \vv{U}_{u}(\vv{x}-\vv{z})\in C$, and by the bounded curvature of $F(\cdot)$, 
\small
\begin{align*}
&F( \vv{z} +  \beta \vv{U}_{u}(\vv{x}-\vv{z})) \le F(\vv{z}) + \beta\partial F(\vv{z})^T\vv{U}_{u}(\vv{x}-\vv{z}) \\
&\hspace{4.5cm}+ \mu_F \beta^2\|\vv{U}_{u}(\vv{x}-\vv{z})\|_2^2\\
&\quad= F(\vv{z}) + \beta\partial F(\vv{z})^T\vv{U}_{u}(\vv{y}-\vv{z})\\
&\qquad+\beta\partial F(\vv{z})^T\vv{U}_{u}(\vv{x}-\vv{y})+ \mu_F\beta^2\|\vv{U}_{u}(\vv{x}-\vv{z})\|_2^2\\
&\quad\le F(\vv{z}) - \beta\epsilon +\beta\partial F(\vv{z})^T\vv{U}_{u}(\vv{x}-\vv{y}) + \mu_F\beta^2\bar{x}_D^2
\end{align*}
\normalsize
where $\bar{x}_D:=\max_{\vv{y},\vv{z}\in C}\|\vv{y}-\vv{z}\|_2$ (the max always exists since $D$ is compact).  By Corollary \ref{cor:one} we can select $\vv{x}\in D$ such that $\partial F(\vv{z})^T\vv{U}_{u}(\vv{x}-\vv{y})\le 0$.  With this choice of $\vv{x}$ it follows that
$
F( \vv{z} + \beta \vv{U}_{u}(\vv{x}-\vv{z}))
\le F(\vv{z}) - \beta\epsilon + \mu_F\beta^2\bar{x}_D^2
$
and the stated result now follows. 

We have the following corollary to Lemma \ref{th:coord} when $F$ is separable:
\begin{corollary}[Descent]\label{cor:sep}
Let $\U\subset 2^{\{1,\cdots,n\}}$ and $F(\vv{z})=\sum_{u\in\U} \vi{F}{u}(\vvu{z}{u})$ with $\vi{F}{u}:\mathbb{R}^{|u|}\rightarrow \mathbb{R}$ convex with bounded curvature with curvature constant $\mu_F$.  Suppose points $\vv{y}$, $\vv{z}\in C=\conv(D)$ exist such that $\vi{F}{u}(\vvu{y}{u})\le \vi{F}{u}(\vvu{z}{u})-\epsilon$ for some $u\in\U$, with $\epsilon>0$, $D$ a $\U$-feasible, compact subset of $\mathbb{R}^n$.  
Selecting $\vv{x} \in\arg\min_{\vv{w}\in D} F(\vvu{z}{u}+ \beta (\vvu{w}{u}-\vvu{z}{u}))$
then $F(\vv{z}+\beta \vv{U}_{u}(\vv{x}-\vv{z}))\le F(\vv{z})-\gamma\beta\epsilon$ provided $\beta\le(1-\gamma)\min\{\epsilon / (\mu_F\bar{x}^2_D),1\}$.   Equivalently, 
\begin{align*}
\vi{F}{u}(\vvu{z}{u}+\beta(\vvu{x}{u}-\vvu{z}{u}))\le \vi{F}{u}(\vvu{z}{u})-\gamma\beta\epsilon
\end{align*}
\end{corollary}
\begin{IEEEproof}
Observe that $F(\vv{z}+\vv{U}_{u}(\vv{y}-\vv{z})) = \sum_{v\in\U\setminus u}\vi{F}{v}( \vvu{z}{v})+\vi{F}{u}(\vvu{y}{u})$ and $F(\vv{z})=\sum_{v\in\U} \vi{F}{v}(\vvu{z}{v})$.  Hence, $F(\vv{z}+\vv{U}_{u}(\vv{y}-\vv{z}))\le F(\vv{z})-\epsilon$ holds when  $\vi{F}{u}(\vvu{y}{u})\le \vi{F}{u}(\vvu{z}{u})-\epsilon$.  Observe also that $F( \vv{z} + \beta \vv{U}_{u}(\vv{x}-\vv{z}))=\sum_{v\in\U\setminus u}\vi{F}{v}( \vvu{z}{v}) + \vi{F}{u}(\vvu{z}{u}+\beta(\vvu{x}{u}-\vvu{z}{u}))$ and $F(\vv{z})-\gamma\beta\epsilon = \sum_{v\in\U}\vi{F}{v}( \vvu{z}{v})-\gamma\beta\epsilon$
\end{IEEEproof}

\begin{lemma}[Unsynchronised Descent]\label{th:decrease}
Let $u_k$, $k=1,2,\dots,N$ be a sequence of elements from admissible update set $\mathcal{U}$ such that $\cup_{k\in\{1,\cdots,N\}}u_k = \{1,\dots,n\}$.   
Let $\{F_k\}$ be a sequence of convex functions with $F_k(\vv{z})=\sum_{u\in\U} \vi{F_k}{u}(\vvu{z}{u})$ where $\vi{F_k}{u}:\mathbb{R}^{|u|}\rightarrow \mathbb{R}$ is convex with bounded curvature and curvature constant $\mu_F$.  Suppose $D$ is $\U$-feasible and compact.
Let $\{\vv{z}_k\}$ be a sequence of vectors with $\vv{z}_{k+1}=\vv{z}_k+ \beta \vv{U}_{u_k} (\vv{x_k}-\vv{z_k})$, $\vv{z}_1\in C=\conv(D)$ and where $\vv{x}_k\in D$ satisfies 
\begin{align*}
&\vi{F}{u_k}_k(\vvu{z}{u_{k}}_k+\beta(\vvu{x}{u_{k}}_k-\vvu{z}{u_{k}}_k))\le 
\vi{F}{u_k}_k(\vvu{z}{u_{k}}_k)+\mu_F\beta^2\bar{x}_D^2
\end{align*}
and, in addition, whenever there exists a point $\vv{y}\in C$ such that $\vi{F}{u_k}_k(\vvu{y}{u_{k}})\le \vi{F}{u_k}_k(\vvu{z}{u_{k}}_k)-\epsilon^\prime$, then
\begin{align*}
&\vi{F}{u_k}_k(\vvu{z}{u_{k}}_k+\beta(\vvu{x}{u_{k}}_k-\vvu{z}{u_{k}}_k))\le \vi{F}{u_k}_k(\vvu{z}{u_{k}}_k)-\gamma\beta\epsilon^\prime
\end{align*}
Suppose that parameter $\beta$ is sufficiently small that
\begin{align*}
\beta\le (1-\gamma)\gamma\min\{\frac{\epsilon^\prime(1+\gamma_1\gamma)}{N\mu_F\bar{x}^2_D},1\}
\end{align*}
with $0<\gamma<1$, $\epsilon^\prime>0$ and that functions $F_k$ change sufficiently slowly that
\begin{align*}
|F_{k+1}(\vv{z}) - F_k(\vv{z})| \le \frac{\gamma_1\gamma}{2N}\beta\epsilon^\prime,\ \forall \vv{z}\in C
\end{align*}
where $\gamma_1\in(0,\gamma/2)$.  Suppose also that a point $\vv{y}\in C$ exists such that $\vi{F_1}{u}(\vvu{y}{u})\le \vi{F_1}{u}(\vvu{z_1}{u})-\epsilon^\prime (1+\gamma_1\gamma\beta)$ for at least one $u\in\mathcal{U}$.
Then 
\begin{align*}
F_{N+1}( \vv{z}_{N+1})<F_1(\vv{z}_1)-\frac{\gamma}{2}\beta\epsilon^\prime
\end{align*}
\end{lemma}
\begin{IEEEproof}
At each step $k\in\{1,\cdots,N\}$ we have two cases to consider.  Namely, either (i) there exists a point $\vv{y}\in C$ such that $\vi{F_k}{u_k}(\vvu{y}{u_{k}})\le \vi{F_k}{u_k}(\vvu{z}{u_{k}})-\epsilon^\prime$ or (ii) there does not.   
In case (i), we have $\vv{x}_k\in D$ such that $\vi{F_k}{u_k}(\vvu{z}{u_{k}}+\beta(\vvu{x}{u_{k}}-\vvu{z}{u_{k}}))\le \vi{F_k}{u_k}(\vvu{z}{u_{k}})-\gamma\beta\epsilon^\prime$.  That is, $F_k( \vv{z}_{k+1})<F_k(\vv{z}_k)-\gamma\beta\epsilon^\prime$.  It follows that $F_{k+1}( \vv{z}_{k+1})<F_k(\vv{z}_k)-(1-\frac{\gamma_1}{N})\gamma\beta\epsilon^\prime$.
In case (ii), $F_k( \vv{z}_{k+1})<F_k(\vv{z}_k)+ \mu_F \beta^2\bar{x}_D^2$.  It follows that $F_{k+1}( \vv{z}_{k+1})<F_k(\vv{z}_k)+ \mu_F \beta^2\bar{x}_D^2+\frac{\gamma_1\gamma}{N}\beta\epsilon^\prime$.
Hence, provided case (i) occurs at least once for $k\in\{1,\cdots,N\}$ then $F_{N+1}(\vv{z}_{N+1})\le F_1(\vv{z}_1)-(1-\gamma_1)\gamma\beta\epsilon^\prime+N\mu_F \beta^2\bar{x}_D^2  \le -\frac{\gamma}{2}\beta\epsilon^\prime$. 

To complete the proof it remains to show that case (i) occurs at least once for $k\in\{1,\cdots,N\}$.  
By assumption there exists a point $\vv{y}\in C$ such that $\vi{F_1}{u}(\vvu{y}{u})\le \vi{F_1}{u}(\vvu{z}{u})-\epsilon^\prime(1+\gamma_1\gamma\beta)$ for at least one $u\in\mathcal{U}$.  Hence, $\vi{F_k}{u}(\vvu{y}{u})\le \vi{F_1}{u}(\vvu{z}{u})-\epsilon^\prime(1+\gamma_1\gamma\beta)+N\frac{\gamma_1\gamma}{N}\beta\epsilon^\prime = \vi{F_1}{u}(\vvu{z}{u})-\epsilon^\prime$ for all $k\in\{1,\cdots,N\}$.
Let $k^*=\min \{k\in \{1,\cdots,N\} : u_{k}\cap u^* = u^*\}$ denote the first time step at which a set $u^*$ satisfying this condition occurs in the sequence.  
We have that $\vvu{z_{k^*}}{u^*}=\vvu{z_1}{u^*}$ since elements $u^*$ of vector $\vv{z}_1$ are not updated at steps $k< k^*$ because $u_k\cap u^*=\emptyset$ when $u_k\ne u^*$ for $u^*$, $u_k\in\bar{\U}$ (note that we do not need to consider $u_{k}\in\U\setminus\bar{\U}$ since $u_{k}\cap u^* = u^*$ for $u_{k}\in\U\setminus\bar{\U}=\{1,\cdots,n\}$ contradicting the assumption that $k< k^*$)  and $k^*$ is the first occurence of $u^*$.  Hence, case (i) occurs at step $k^*$ and we are done.
\end{IEEEproof}

\begin{lemma}[Descent into $C_{(2p+1)\epsilon}$ Ball]\label{th:intoball}
Let $u_k$, $k=1,2,\cdots$ be a sequence of elements from admissible update set $\mathcal{U}$.   Suppose we can partition this sequence into disjoint intervals $\{k_i,\cdots,k_{i+1}\}$, $i=1,2,\dots$ with $k_i> k_{i-1}$
and such that over each interval $\cup_{k\in\{k_j,k_{j+1}\}} u_k = \{1,\cdots,n\}$.  Suppose the largest size of interval $N=\max_{j}(k_{j+1}-k_j+1)$ is finite and that  
the conditions in Lemma \ref{th:decrease} hold over each interval with $\epsilon^\prime=\epsilon/(1+\gamma_1\gamma\beta)$.  
Then for $k_i$ sufficiently large we have that, 
\begin{align*}
F_{k_{i}}(\vv{z}_{k_i})-F_{k_{i}}(\vv{y}_{k_{i}}^*)\le \left({|\U|}+ 1\right)\epsilon
\end{align*}
where $\vv{y}^*_k\in\arg\min_{\vv{z}\in C} F_k(\vv{z})$.
\end{lemma}
\begin{IEEEproof}
Consider interval $\{k_j,\cdots,k_{j+1}\}$.  We proceed by considering two cases: 
Case (i): $\exists u\in\mathcal{U}$ such that $\vi{F_{k_j}}{u}(\vvu{z_{k_j}}{u}) - \vi{F_{k_j}}{u}(\vvu{y^*_{k_j}}{u}) \ge \epsilon$.   Applying  Lemma \ref{th:decrease} with parameter $\epsilon^\prime=\epsilon/(1+\gamma_1\gamma\beta)$, we have that $F_{k_{j+1}}(\vv{z}_{k_{j+1}}) - F_{k_j}(\vv{z}_{k_j}) \le - \frac{\gamma}{2}\beta \epsilon^\prime <0$.   That is, $F_{k_{j+1}}(\cdot)$ decreases monotonically.
Case (ii):  $\vi{F_{k_j}}{u}(\vvu{z_{k_j}}{u}) - \vi{F_{k_j}}{u}(\vvu{y^*_{k_j}}{u}) < \epsilon$ $\forall u\in\mathcal{U}$.  It follows that $F_{k_j}(\vv{z}_{k_j}) < F_{k_j}(\vv{y}_{k_j}^*)+{|\U|}\epsilon$.  Applying the same argument as for case (ii) of the proof of Lemma \ref{th:decrease}, it follows that $F_{k_{j+1}}(\vv{z}_{k_{j+1}}) < F_{k_{j}}(\vv{y}_{k_j}^*)+{|\U|}\epsilon + N\mu_F\beta^2\bar{x}_D^2 +\gamma_1\gamma\beta\epsilon$.  
Further, $F_{k_{j}}(\vv{y}_{k_{j}}^*) \le F_{k_{j}}(\vv{y}_{k_{j+1}}^*) \le F_{k_{j+1}}(\vv{y}_{k_{j+1}}^*)+\gamma_1\gamma\beta\epsilon$.  Hence, $F_{k_{j+1}}(\vv{z}_{k_{j+1}}) < F_{k_{j+1}}(\vv{y}_{k_{j+1}}^*)+{|\U|}\epsilon^\prime + N\mu_F\beta^2\bar{x}_D^2 +2\gamma_1\gamma\beta\epsilon$
Using the stated choice of $\beta$ yields
\begin{align*}
F_{k_{j+1}}(\vv{z}_{k_{j+1}}) &\le F_{k_{j+1}}(\vv{y}_{k+1}^*)+{|\U|}\epsilon 
+ (1-\gamma+2\gamma_1)\gamma\beta\epsilon\\
&\le F_{k_{j+1}}(\vv{y}_{k+1}^*)+{|\U|}\epsilon + \gamma\beta\epsilon
\end{align*}
Observe that when $F_{k_j}(\vv{z}_{k_j}) - F_{k_j}(\vv{y}_{k_j}^*)\ge {|\U|}\epsilon$ then there must exist 
a $u\in\mathcal{U}$ such that $\vi{F_{k_j}}{u}(\vvu{z_{k_j}}{u}) - \vi{F_{k_j}}{u}(\vvu{y^*_{k_j}}{u}) \ge \epsilon$.  We therefore have that $F_{k+1}(\vv{z}_{k+1})$ is strictly decreasing when $F_{k_j}(\vv{z}_{k_j})-F_{k_j}(\vv{y}_{k_j}^*) \ge {|\U|}\epsilon$ and otherwise any increase is uniformly upper bounded by ${|\U|}\epsilon + \gamma\beta\epsilon$.  It follows that for all $k$ sufficiently large $F_{k_{i}}(\vv{z}_{k_i})-F_{k_{i}}(\vv{y}_{k_{i}}^*) \le {|\U|}\epsilon + \gamma\beta\epsilon\le\left({|\U|}+ 1\right)\epsilon$ as claimed.
\end{IEEEproof}

\subsection*{Proof of Theorem \ref{th:unsynch}}

By the bounded curvature of $\vi{F_k}{u_k}$ and Corollary \ref{cor:one},
\small
\begin{align*}
&\vi{F_k}{u_k}( \vvu{z}{u_{k}} + \beta (\vvu{x}{u_{k}}-\vvu{z}{u_{k}})) \\
&\qquad\le \vi{F_k}{u_k}(\vvu{z}{u_{k}}) + \beta\partial \vi{F_k}{u_k}(\vvu{z}{u_{k}})^T (\vvu{x}{u_{k}}-\vvu{z}{u_{k}}) \\
&\qquad\qquad+ \mu_F \beta^2\|\vvu{x}{u_{k}}-\vvu{z}{u_{k}}\|_2^2\\
&\qquad\stackrel{(a)}{\le} \vi{F_k}{u_k}(\vvu{z}{u_{k}}) + \mu_F \beta^2\|\vvu{x}{u_k}-\vvu{z}{u_k}\|_2^2 \\
&\qquad\le \vi{F_k}{u_k}(\vvu{z}{u_{k}}) + \mu_F \beta^2\bar{x}_D^2
\end{align*}
\normalsize
where $(a)$ is from the minimality of the choice of $\vv{x}_k$.   Combining with Corollary \ref{cor:sep}, the conditions of Lemma \ref{th:decrease} are satisfied with $\epsilon^\prime=\epsilon/(1+\gamma_1\gamma\beta)$ for the stated choice of $\vv{x}_k$, $k=1,2,\cdots$.  Theorem \ref{th:unsynch} follows by application of Lemma \ref{th:intoball}.

\begin{lemma}[Descent]\label{th:coord2}
Let $F(\vv{z})$ be convex with bounded curvature on $C$ with curvature constant $\mu_F$.  Suppose points $\vv{y}$, $\vv{z}\in C=\conv(D)$, $D$ a $\U$-feasible compact subset of $\mathbb{R}^n$ exist such that $F(\vv{z}+\vv{U}_{u}(\vv{y}-\vv{z}))\le F(\vv{z})-\epsilon$, $\epsilon>0$ and $u\in\mathcal{U}\subset 2^{\{1,\cdots,n\}}$.  
Selecting
\begin{align*}
\vv{x} &\in\arg\min_{\vv{w}\in D} \partial F(\vv{z})^T\vv{U}_{u}\vv{w}
\end{align*}
then
\begin{align*}
F( \vv{z} + \beta \vv{U}_{u}(\vv{x}-\vv{z}))<F(\vv{z})-\gamma\beta\epsilon
\end{align*}
 provided $\beta\le (1-\gamma)\min\{\epsilon / (\mu_F\bar{x}^2_D),1\}$ $0<\gamma<1$ where $\bar{x}_D:=\max_{\vv{y},\vv{z}\in \conv(D)}\|\vv{y}-\vv{z}\|_2$.
\end{lemma}
\begin{IEEEproof}
Since $D$ consists of the extreme points of set $C$ observe that $\arg\min_{\vv{x}\in C} F_k(\vv{z}_k)+\partial F_k(\vv{z}_k)^T\vv{U}_{u}(\vv{x}-\vv{z}_k) = \arg\min_{\vv{x}\in C}\partial F_k(\vv{z}_k)^T\vv{U}_{u}(\vv{x}-\vv{z}_k)  = \arg\min_{\vv{x}\in D}\partial F_k(\vv{z}_k)^T\vv{U}_{u}(\vv{x}-\vv{z}_k)$ (since the solution to a linear programme lies at an extreme point).   
By convexity and the minimality of $\vv{x}$,
$F(\vv{z}) +  \partial F(\vv{z})^T\vv{U}_{u}(\vv{x}-\vv{z})
\le F(\vv{z}) +  \partial F(\vv{z})^T\vv{U}_{u}(\vv{y}-\vv{z})
\le F(\vv{z})-\epsilon$.
Hence, $\partial F(\vv{z})^T\vv{U}_{u}(\vv{x}-\vv{z}) \le -\epsilon$.    Since $D$ is $\U$-feasible, for any $\vv{x}\in D$ we have $ \vv{z} + \beta \vv{U}_{u}(\vv{x}-\vv{z})\in C$, and by the bounded curvature of $F$, 
$F( \vv{z} +  \beta \vv{U}_{u}(\vv{x}-\vv{z})) 
\le F(\vv{z}) + \beta\partial F(\vv{z})^T\vv{U}_{u}(\vv{x}-\vv{z}) 
+ \mu_F \beta^2\|\vv{U}_{u}(\vv{x}-\vv{z})\|_2^2
\le F(\vv{z}) - \beta\epsilon + \mu_F\beta^2\bar{x}_D^2 \le F(\vv{z})  -\gamma\beta\epsilon
$.
\end{IEEEproof}

\subsection*{Proof of Theorem \ref{th:FWunsynch}}

By Corollary \ref{cor:one} and Lemma \ref{th:coord2} the conditions of Lemma \ref{th:decrease} are satisfied for the stated choice of $\vv{x}_k$, $k=1,2,\cdots$.  Theorem \ref{th:FWunsynch} now follows by application of Lemma \ref{th:intoball}.

\subsection{Approximate Solutions}

\begin{lemma}[Bounded Multipliers]\label{th:bmulti} 
Consider the setup of problem $P$ and update $
\vv{\lambda}_{k+1} = [\vv{\lambda}_k + \alpha \vv{g}(\vv{z}_k)]^+$ where $\vv{z}_k \in  C_\epsilon (\vv{\mu}_k) := \left\{ \vv{z} \in C \mid L(\vv{z}, \vv{\mu}_k) - q(\vv{\mu}_k) \le \epsilon \right\}$, $\epsilon \ge 0$ and $\| \vv{\lambda}_k - \vv{\mu}_k \|_\infty \le \alpha \sigma_0$ for all $k$ with $\alpha > 0$, $\sigma_0 \ge 0$. Suppose the Slater condition is satisfied and that $\vv{\lambda}_1 \in \mathbb R^m_+$. Then, for all $k=1,2,\dots$ we have that
\begin{align}
&  \|\vv{\lambda}_{k}\|_2 \le \bar \lambda:= 2 \mathcal Q  + \max \{ \| \vv{\lambda}_{1}\|_2 , \mathcal Q  + \alpha m {\bar g} \}
\end{align}where ${\bar g}:= \max_{\vv{z} \in C} \| \vv{g}(\vv{z})\|_\infty$ and $\mathcal Q$ is given in Lemma \ref{th:setq} with $\delta:= \alpha m^2 ({\bar g}^2/2 + \sigma_0 {\bar g} ) + \epsilon $.
\end{lemma}
%
\begin{IEEEproof}
First let $\vv{\theta} \in \mathbb R^m$ and see that
\small
\begin{align}
&\| \vv{\lambda}_{k+1}-\vv{\theta}\|_2^2 
 \le \|[\vv{\lambda}_k + \alpha \vv{g}(\vv{z}_{k})]^{+}-\vv{\theta}\|_2^2 \label{eq:low1} \\
&\quad  \le \| \vv{\lambda}_k + \alpha \vv{g}(\vv{z}_{k})-\vv{\theta}\|_2^2  \notag \\
&\quad = \| \vv{\lambda}_k-\vv{\theta}\|_2^2 + 2\alpha(\vv{\lambda}_k-\vv{\theta})^T\vv{g}(\vv{z}_{k}) + \alpha^2 \| \vv{g}(\vv{z}_{k}) \|_2^2 \notag \\
&\quad \le \| \vv{\lambda}_k-\vv{\theta}\|_2^2 + 2\alpha(\vv{\lambda}_k-\vv{\theta})^T\vv{g}(\vv{z}_{k}) + \alpha^2 m {\bar g}^2 \label{eq:low2} \\
&\quad = \| \vv{\lambda}_k-\vv{\theta}\|_2^2 + 2\alpha(L(\vv{z}_k,\vv{\lambda}_k)-L(\vv{z}_k,\vv{\theta})) + \alpha^2 m {\bar g}^2 \label{eq:low3}
\end{align} 
\normalsize
where $(\ref{eq:low1})$ follows since $\vv{\theta} \succeq \vv{0}$, $(\ref{eq:low2})$ from the fact that $\|\vv{g}(\vv{z})\|^2_2 \le m{\bar g}^2$ for all $\vv{z} \in C$, (\ref{eq:low3}) since $(\vv{\lambda}_k - \vv{\theta})^T \vv{g}(\vv{z}_k) = f(\vv{z}_k) - f(\vv{z}_k) + (\vv{\lambda}_k - \vv{\theta})^T \vv{g}(\vv{z}_k) = L(\vv{z}_k,\vv{\lambda}_k)-L(\vv{z}_k,\vv{\theta})$. Further, see that since $\| \vv{\lambda}_k - \vv{\mu}_k \|_\infty \le \alpha \sigma_0$ for all $k$ we have that  $| L(\vv{z}_k, \vv{\lambda}_k) - L(\vv{z}_k , \vv{\mu}_k) | = | (\vv{\lambda}_k - \vv{\mu}_k)^T \vv{g}(\vv{z}_k)| \le \| \vv{\lambda}_k  - \vv{\mu}_k \|_2 \| \vv{g}(\vv{z}_k) \|_2 \le \alpha m \sigma_0 \bar g$. Using the latter bound in (\ref{eq:low3}) and rearranging terms yields
\begin{align}
 &\| \vv{\lambda}_{k+1} - \vv{\theta} \|_2^2  - \| \vv{\lambda}_k - \vv{\theta} \|_2^2 \notag\\ 
& \qquad \le  \alpha^2 m (\bar g^2 + 2\sigma_0 \bar g)  + 2 \alpha ( L(\vv{z}_k, {\vv{\mu}}_k) - L(\vv{z}_k,\vv{\theta})) \label{eq:stepstoslackness} \\
&  \qquad  =  \alpha^2 m ({\bar g}^2 + 2\sigma_0 \bar g)  + 2 \alpha ( q({\vv{\mu}}_k) + \epsilon - L(\vv{z}_k,\vv{\theta}))	\label{eq:steptotheorem}
\end{align}
where the last equation follows from the fact that $\vv{z}_k \in C_\epsilon (\vv{\mu}_k)$  {i.e.}, $L(\vv{z}_k, {\vv{\mu}}_k) \le q({\vv{\mu}}_k) + \epsilon$. Letting $\vv{\theta} = \vv{\lambda}^\star$ and using the fact that $L(\vv{z}_k,\vv{\lambda}^\star) \ge \min_{\vv{z} \in C} L(\vv{z},{\vv{\lambda}}_k) = q(\vv{\lambda}^\star)$ we obtain
\begin{align}
 &\| \vv{\lambda}_{k+1} - \vv{\lambda}^\star \|_2^2  - \| \vv{\lambda}_k - \vv{\lambda}^\star \|_2^2 \notag \\
 &\qquad\le  \alpha^2 m ({\bar g}^2 + 2\sigma_0 \bar g)  + 2 \alpha ( q({\vv{\mu}}_k) + \epsilon - q(\vv{\lambda}^\star))	\label{eq:diffslackboundfirst}
\end{align}

Now let $Q_\delta : = \{ {\vv{\mu}} \succeq \vv{0} : q({\vv{\mu}}) \ge q(\vv{\lambda}^\star) - \delta )\}$ and consider two cases. Case (i) $({\vv{\mu}}_k \notin Q_{\delta})$. Then $q({\vv{\mu}}_k) - q(\vv{\lambda}^\star) < - \delta$ and from (\ref{eq:diffslackboundfirst}) we have that $\| \vv{\lambda}_{k+1} - \vv{\lambda}^\star \|_2^2 < \| \vv{\lambda}_{k} - \vv{\lambda}^\star \|_2^2$, \emph{i.e.,} 
 \begin{align*}
 \| \vv{\lambda}_{k+1} - \vv{\lambda}^\star \|_2  - \| \vv{\lambda}_k - \vv{\lambda}^\star \|_2 < 0
 \end{align*} and so $\vv{\lambda}_k$ converges into a ball around $\vv{\lambda}^\star$ when ${\vv{\mu}}_k \in Q_\delta$.  Case (ii) $({\vv{\mu}}_k \in Q_{\delta})$. See that  $\| \vv{\lambda}_{k+1} - \vv{\lambda}^\star \| = \| [\vv{\lambda}_k + \alpha \vv{g}(\vv{z}_{k})]^+ - \vv{\lambda}^\star \|_2 \le \| \vv{\lambda}_k + \alpha \vv{g}(\vv{z}_{k}) - \vv{\lambda}^\star \|_2 \le \| \vv{\lambda}_k \|_2 + \| \vv{\lambda}^\star \|_2 + \alpha m {\bar g}$. Next recall that when the Slater condition holds by Lemma \ref{th:setq} we have that $\| \vv{\lambda} \|_2 \le \mathcal Q$ for all $\vv{\lambda} \in Q_\delta$. Therefore, 
 \begin{align*}
 \| \vv{\lambda}_{k+1} - \vv{\lambda}^\star \|_2 \le \mathcal Q + \alpha m {\bar g}.
 \end{align*}Combining both cases we have that
\begin{align*}
 \| \vv{\lambda}_{k+1} - \vv{\lambda}^\star \|_2 &   \le \max \{ \| \vv{\lambda}_{1} - \vv{\lambda}^\star \|_2 ,2 \mathcal Q + \alpha m {\bar g} \} \\
&   \le \max \{ \| \vv{\lambda}_{1}\|_2 + \| \vv{\lambda}^\star \|_2 , 2 \mathcal Q + \alpha m {\bar g} \} \\
&   \le \mathcal Q+ \max \{ \| \vv{\lambda}_{1}\|_2 ,  \mathcal Q + \alpha m {\bar g} \}.
\end{align*}Finally, since we have that  $\|\vv{\lambda}_{k+1} - \vv{\lambda}^\star \|_2 \ge | \| \vv{\lambda}_{k+1}\|_2 - \| \vv{\lambda}^\star\|_2 | \ge   \| \vv{\lambda}_{k+1}\|_2 - \| \vv{\lambda}^\star\|_2 \ge  \|\vv{\lambda}_{k+1} \|_2 - \mathcal Q$ the stated result now follows.
\end{IEEEproof}

\subsection*{Proof of Theorem \ref{th:maintheorem}}

Applying the bound in (\ref{eq:stepstoslackness}) recursively from $i=1,\dots,k$ we have that 
$\sum_{i=1}^k  \| \vv{\lambda}_{i+1} - \vv{\theta} \|_2^2  - \|\vv{\lambda}_i - \vv{\theta} \|_2^2  =   \| \vv{\lambda}_{k+1} - \vv{\theta} \|_2^2  - \| \vv{\lambda}_1 - \vv{\theta} \|_2^2 
  \le  \alpha^2 m ({\bar g}^2 + \sigma_0 \bar g) k  + 2 \alpha \sum_{i=1}^k ( L(\vv{z}_i,{\vv{\mu}}_i) - L(\vv{z}_i,\vv{\theta}))$. Rearranging terms  and dividing by $2 \alpha k$ yields 
  \small
\begin{align}
& - \frac{\| \vv{\lambda}_1 - \vv{\theta}\|_2^2}{2\alpha k } - {\alpha} m (\frac{{\bar g}^2}{ 2} + \sigma_0 \bar g) 
 \le \frac{1}{k} \sum_{i=1}^{k} (L(\vv{z}_i,\vv{\mu}_i) -L(\vv{z}_i,\vv{\theta})). \label{eq:subgradup}
\end{align}
\normalsize

We now consider two cases. Case (i). Let $\vv{\theta} = \vv{\lambda}^\star$ and see that by the saddle-point property we have that $L(\vv{z},\vv{\theta}) = L(\vv{z},\vv{\lambda}^\star) \ge \min_{\vv{z} \in C} L(\vv{z},\vv{\lambda}^\star) =  q(\vv{\lambda}^\star) = f^\star$. Further, let $ \vv{\lambda}^\diamond_k  := \frac{1}{k} \sum_{i=1}^k \vv{\lambda}_i$  and see $ \frac{1}{k} \sum_{i=1}^k L(\vv{z}_i, \vv{\mu}_i) \le \frac{1}{k} \sum_{i=1}^k q(\vv{\mu}_i) + \epsilon \le q(\vv{\mu}^\diamond_k) + \epsilon$ by the concavity of $q$. Combining both bounds yields
\begin{align}
 - \frac{m \bar \lambda^2 }{\alpha k} -  {\alpha} m({\bar g}^2/2 + \sigma_0 \bar g) - \epsilon \le q({\vv{\mu}}^{\diamond}_k) -  q(\vv{\lambda}^\star) \le 0 \label{eq:qavgparteq}
\end{align}
Case (ii). Let $\vv{\theta} = \vv{\mu}^\diamond_k$ and see that by the convexity of $L(\cdot, \vv{\theta})$ we have that $\frac{1}{k} \sum_{i=1}^k L(\vv{z}_i , \vv{\theta}) = \frac{1}{k} \sum_{i=1}^k L(\vv{z}_i , {\vv{\mu}}^\diamond_k) \ge  L(\vv{z}^\diamond_k, {\vv{\mu}}^\diamond_k)$. Using again the fact that $ \frac{1}{k} \sum_{i=1}^k L(\vv{z}_i, \vv{\mu}_i) \le q(\vv{\mu}^\diamond_k) + \epsilon$ follows that
\begin{align*}
& - \frac{m \bar \lambda^2}{\alpha k} -  {\alpha} m (\bar g^2/2 + \sigma_0 \bar g) - \epsilon \le q({\vv{\mu}}^\diamond_k)  - L(\vv{z}^\diamond_{k},{\vv{\mu}}^\diamond_k) \stackrel{(a)}{\le} 0
\end{align*}where $(a)$ holds since $L(\vv{z}^\diamond_{k},{\vv{\mu}}^\diamond_k) \ge \min_{\vv{z}\in C} L(\vv{z},{\vv{\mu}}^\diamond_k) := q({\vv{\mu}}^\diamond_k)$. Multiplying the latter bound by $-1$ and combining it with (\ref{eq:qavgparteq}) yields 
\begin{align}
| L(\vv{z}^\diamond_{k},{\vv{\mu}}^\diamond_k) - f^* | \le \epsilon +  \alpha m (\bar g^2/ 2 + \sigma_0 \bar g) + \frac{m \bar \lambda^2}{\alpha k} \label{eq:lagrangianbound}
\end{align}

We now proceed the upper and lower bound $L(\vv{z}^\diamond_k, \vv{\lambda}^\diamond_k) - f(\vv{z}^\diamond_k) =  (\vv{\lambda}^\diamond_k)^T \vv{g}(\vv{z}^\diamond_k)$. For the lower bound let $\vv{\theta} = \vv{0}$ in (\ref{eq:subgradup}) and observe that 
\begin{align}
& - \frac{\| \vv{\lambda}_1\|_2^2}{2\alpha k } - {\alpha} m ({\bar g}^2/2 + \sigma_0 \bar g)   \le \frac{1}{k} \sum_{i=1}^{k} \vv{\mu}_i^T \vv{g}(\vv{z}_i) 
\end{align}The LHS of the last equation does not depend on $\vv{z}_i$, {i.e.}, we can fix $\vv{z}_i = \vv{z}^\diamond_k$ for all $i=1,\dots,k$, and so 
 \begin{align*}
- \frac{\| \vv{\lambda}_{1} \|^2_2}{2 \alpha k } - {\alpha} m {\bar g}^2  & \le  \vv{g}(\vv{z}^\diamond_k)  \sum_{i=1}^k  {\vv{\lambda}}_i  =  ({\vv{\lambda}}^\diamond_k)^T \vv{g}(\vv{z}^\diamond_k).  
 \end{align*}
To upper bound $({\vv{\lambda}}^\diamond_k)^T \vv{g}(\vv{z}^\diamond_k)$ see that
\small
\begin{align}
 \vv{\lambda}_{k+1} & = \alpha  \max \left\{ \max_{1 \le j \le k} \sum_{i=j}^k \vv{g}(\vv{z}_{i}), \left[ \sum_{i=1}^k \vv{g}(\vv{z}_i)  + \vv{\lambda}_1 \right]^+  \right\} \notag \\
& \ge   \alpha \max \left\{  \sum_{i=1}^k \vv{g}(\vv{z}_{i}), \left[\sum_{i=1}^k \vv{g}(\vv{z}_i)  + \vv{\lambda}_1 \right]^+  \right\} \label{eq:sumavgslack} \\
&  \ge \alpha \sum_{i=1}^k \vv{g}(\vv{z}_i)
  \ge \alpha k  \vv{g}(\vv{z}^\diamond_k) \label{eq:convavgslack}
\end{align}
\normalsize
where (\ref{eq:sumavgslack}) follows from the fact that $\sum_{i=1}^k \vv{g}(\vv{z}_i) \le \max_{1 \le j \le k } \sum_j^k \vv{g}(\vv{z}_i)$ and (\ref{eq:convavgslack}) follows from the convexity of $\vv{g}(\cdot)$. Multiplying both side by $\vv{\lambda}^\diamond_k$, dividing by $\alpha k$ and using the fact that $\vv{\lambda}^\diamond_k = \frac{1}{k} \sum_{i=1}^k \vv{\lambda}_i \preceq \bar \lambda \vv{1}$ yields
\begin{align}
- \frac{m \bar \lambda^2}{2\alpha k } - {\alpha} m {\bar g}^2 / 2  \le  ({\vv{\lambda}}^\diamond_k)^T \vv{g}(\vv{z}^\diamond_k) \le \frac{m \bar \lambda^2}{\alpha k} \label{eq:slackbound}
\end{align}Combining (\ref{eq:slackbound}) with (\ref{eq:lagrangianbound}) and using the fact that $\| \vv{\lambda}^\diamond_k - \vv{\mu}^\diamond_k \|_\infty = \| \frac{1}{k} \sum_{i=1}^k (\vv{\lambda}_i - \vv{\mu}_i) \|_\infty \le \alpha \sigma_0$ yields the stated result.


\subsection{Discrete Actions: Proof of Theorem \ref{th:discretesequence}}

We begin by noting that since $\1^T\vv{a}_k=1=\1^T\vv{b}_k$ then $\1^T \vv{S}_k=0$ for all $k=1,2,\cdots$.   Also note that since $\vv{a}_k\in\Delta$ all elements of $\vv{a}_k$ are non-negative and at least one element must be non-zero since $\1^T\vv{a}_k=1$.

We now proceed by induction to show that there always exists a choice of $\vv{b}_{k+1}$ such that $\vv{S}_k\succeq -\1$, $k=1,2,\cdots$.  
%
When $k=1$ let element $j$ of $\vv{a}_1$ be positive (as already noted, at least one such element exists).   Selecting $\vv{b}_1=\vv{e}_j$ then it follows that $-1< \vi{a}{j}_1-\vi{b}{j}_1\le 0$ and so $-1\prec \vv{a}_1-\vv{b}_1 \prec 1$.  That is, $\vv{S}_1\succeq -\1$.   
Suppose now that $\vv{S}_k\succ -\1$.  We need to show that $\vv{S}_{k+1}\succeq -\1$.   Now $\vv{S}_{k+1}=\vv{S}_k+\vv{a}_{k+1}-\vv{b}_{k+1}$.   Since $\vv{S}_k\succeq -\1$, $\vi{S}{j}_k \ge -1$ $\forall j=1,\cdots,|D|$.   Also, $\1^T \vv{S}_k = 0$, so either all elements are 0 or at least one element is positive.   If they are all zero then we are done (we are back to the $k=1$ case).    
Otherwise, since all elements of $\vv{a}_{k+1}$ are non-negative then at least one element of $\vv{S}_k+\vv{a}_{k+1}$ is positive.  Let element $\vi{S}{j}_k+\vi{a}{j}_{k+1}$ be the largest positive element of $\vv{S}_k+\vv{a}_{k+1}$.   Selecting $\vv{b}_{k+1}=\vv{e}_j$ then it follows that $\vi{S}{j}_k+\vi{a}{j}_{k+1}+\vi{b}{j}_{k+1}\ge -1$.   
All elements of $\vv{S}_k+\vv{a}_{k+1}$ are therefore lower bounded by $-1$ since  $\vv{S}_k\succeq -\1$ and the elements of $\vv{a}_{k+1}$ are non-negative as are the elements $\vi{b}{i}_{k+1}=\vi{e}{i}_j$, $i\ne j$.  That is, $\vv{S}_{k+1}\succeq -\1$.

We now have that $\vv{b}_{k+1}$ can always be selected such that $\vv{S}_k\succeq -\underline{S}\1$ provided $\underline{S}\ge 1$, and also that $\1^T \vv{S}_k=0$.  Since $\1^T \vv{S}_k=0$ either $\vv{S}$ is zero or at least one element is positive.   Since $\vv{S}_k\succeq -\underline{S}\1$ and at most $|D|-1$ elements are negative then the sum over the negative elements is lower bounded by $-(|D|-1)\underline{S}$.  Since $\1^T \vv{S}_k=0$ it follows that the sum over the positive elements must be upper bounded by $|D|-1$.  Hence, $\| \vv{S}_k\|_\infty \le (|D|-1)\underline{S}$.
Finally, observe that $\|\sum_{i=1}^k \vv{z}_i-\vv{x}_i\|_\infty = \|\vv{X} \vv{S}_k\|_\infty \le  \|\vv{X}\|_\infty\|\vv{S}_k\|_\infty \le (|D|-1)\underline{S}\|\vv{X}\|_\infty$ and we are done.

\bibliographystyle{unsrt}
\bibliography{bib}

%
%

\end{document}